\documentclass[10pt]{article}
\usepackage{amsmath, amssymb}
\usepackage{amsthm} 
\usepackage[all]{xy}

\theoremstyle{plain} 

\theoremstyle{definition}

\def\Sch{\text{\rm Sch}}

\def\bgn{\begin}

\def\CL{\text{\rm CL}}

\def\GL{\text{\rm GL}}

\def\J{{\mathcal J}}

\def\1{{[1]}}
\def\2{{[2]}}
\def\3{{[3]}}
\def\({\left(}
\def\){\right)}
\def\s-circ{\,{\scriptstyle{\circ}}\,}
\def\<<{<\negthinspace \negthinspace<}
\def\Ad{\text{\rm Ad}}

\def\even{\text{\rm even}}
\def\bgn{\begin}

\def\endaln{\end{align}}

\def\<{<\negthinspace \negthinspace <}

\def\({\left(}
\def\){\right)}
\def\Im{\text{\rm Im}}
\def\Re{\text{\rm Re}}
\def\[{\big[\neg\big[}
\def\]{\big]\neg\big]}
\def\al{\al}

\def\M{{\mathcal M}}

\def\tr{\text{\rm tr}}
\def\a{\alpha}

\def\e{\varepsilon}
\def\gam{\gamma}
\def\Gam{\Gamma}

\def\del{\delta}
\def\lam{\lambda}
\def\ome{\omega}
\def\Ome{\Omega}
\def\sig{\sigma}

\def\A{{\mathcal A} }

\def\Diff{\text{\rm Diff}_0}

\def\R{\Bbb R}
\def\C{\Bbb C}

\def\M{\frak M}

\def\w{\wedge}
\def\({\left(}
\def\){\right)}

\def\neg{\negthinspace}
\def\h{\hat}

\def\til{\tilde}
\def\wtil{\widetilde}

\def\ol{\overline}
\def\pa{\partial}

\def\ran{\rangle} 
\def\lan{\langle}

\def\arrow{\longrightarrow}

\def\:{\, :\,}
\def\CL{\text{\rm CL}}
\def\TT{T\oplus T^*}
\def\complex{generalized complex }
\def\K\"ahler{generalized K\"ahler}

\def\10{\displaystyle L^{10}}
\def\2{\displaystyle L^2}
\def\c0{\displaystyle C^0}
\def\dstyle{\displaystyle}
\def\10{\displaystyle L^{10}}
\def\2{\displaystyle L^2}
\def\del{\delta}
\def\del2{\displaystyle L^2_{0,\delta}}
\def\c0{\displaystyle C^0}
\def\dstyle{\displaystyle}

\def\del{\delta}

\def\K{{\mathcal K}}
\def\M-A{\text{\rm Monge-Amp\`ere}}

\def\M-A{\text{\rm Monge-Amp\`ere}}
\def\[{\big[\,}
\def\]{\,\big]}

\def\id{\text{\rm id}}

\def\tt{\scriptscriptstyle T\oplus T^*}
\def\GR{\text{\rm GR}}

\makeatletter
\def\even{\text{\rm even}}
\def\odd{\text{\rm odd}}
\def\GK{\text{\rm generalized K\"ahler }}

\def\TT{T_M\oplus T^*_M}

\def\M{{\mathcal M}}
\def\Diff{\text{\rm Diff}}
\def\ol{\overline}
\def\part{\partial}
\def\Ham{\text{\rm Ham}}
\def\tt{\scriptstyle  T\oplus T*}

\def\Sch{\text{\rm Sch}}

\def\bgn{\begin}

\def\CL{\text{\rm CL}}

\def\GL{\text{\rm GL}}

\def\J{{\mathcal J}}

\def\1{{[1]}}
\def\2{{[2]}}
\def\3{{[3]}}
\def\({\left(}
\def\){\right)}
\def\s-circ{\,{\scriptstyle{\circ}}\,}
\def\<<{<\negthinspace \negthinspace<}
\def\Ad{\text{\rm Ad}}

\def\even{\text{\rm even}}
\def\bgn{\begin}

\def\endaln{\end{align}}

\def\<{<\negthinspace \negthinspace <}

\def\({\left(}
\def\){\right)}
\def\Im{\text{\rm Im}}
\def\Re{\text{\rm Re}}
\def\[{\big[\neg\big[}
\def\]{\big]\neg\big]}
\def\al{\al}

\def\M{{\mathcal M}}

\def\tr{\text{\rm tr}}
\def\a{\alpha}

\def\e{\varepsilon}
\def\gam{\gamma}
\def\Gam{\Gamma}

\def\del{\delta}
\def\lam{\lambda}
\def\ome{\omega}
\def\Ome{\Omega}
\def\sig{\sigma}

\def\A{{\mathcal A} }

\def\Diff{\text{\rm Diff}_0}

\def\R{\Bbb R}
\def\C{\Bbb C}

\def\M{\frak M}

\def\w{\wedge}
\def\({\left(}
\def\){\right)}

\def\neg{\negthinspace}
\def\h{\hat}

\def\til{\tilde}
\def\wtil{\widetilde}

\def\ol{\overline}
\def\pa{\partial}

\def\ran{\rangle} 
\def\lan{\langle}

\def\arrow{\longrightarrow}

\def\:{\, :\,}
\def\CL{\text{\rm CL}}
\def\TT{T\oplus T^*}
\def\complex{generalized complex }
\def\K\"ahler{generalized K\"ahler}

\def\10{\displaystyle L^{10}}
\def\2{\displaystyle L^2}
\def\c0{\displaystyle C^0}
\def\dstyle{\displaystyle}
\def\10{\displaystyle L^{10}}
\def\2{\displaystyle L^2}
\def\del{\delta}
\def\del2{\displaystyle L^2_{0,\delta}}
\def\c0{\displaystyle C^0}
\def\dstyle{\displaystyle}

\def\del{\delta}

\def\K{{\mathcal K}}
\def\M-A{\text{\rm Monge-Amp\`ere}}

\def\M-A{\text{\rm Monge-Amp\`ere}}
\def\[{\big[\,}
\def\]{\,\big]}

\def\id{\text{\rm id}}

\def\even{\text{\rm even}}
\def\odd{\text{\rm odd}}
\def\GK{\text{\rm generalized K\"ahler }}

\def\TT{T_M\oplus T^*_M}

\def\M{{\mathcal M}}
\def\Diff{\text{\rm Diff}}
\def\ol{\overline}
\def\part{\partial}
\def\Ham{\text{\rm Ham}}

\title{Scalar curvature as moment map \\in generalized K\"ahler geometry
} 
\author{Ryushi Goto}

\date{} 
\begin{document}
\maketitle

\begin{abstract}
It is known that the scalar curvature arises as the moment map in K\"ahler geometry. In pursuit of the analogy,  we develop the moment map framework in generalized K\"ahler geometry of symplectic type. Then we establish  the definition of the scalar curvature on a generalized K\"ahler manifold of symplectic type from the moment map view point.  We also obtain the generalized Ricci form which is a representative of the first Chern class of the anticanonical line bundle.  
We show that infinitesimal deformations of generalized K\"ahler structures with constant generalized scalar curvature are finite dimensional on a compact manifold.
Explicit descriptions of the generalized Ricci form and the generalized scalar curvature are given on a generalized K\"ahler manifold of type $(0,0)$.
Poisson structures constructed from a K\"ahler action of $T^m$ on a K\"ahler-Einstein manifold give
rise to intriguing deformations of generalized K\"ahler-Einstein structures.
In particular, the anticanonical divisor of three lines on $\C P^2$ in general position yields nontrivial examples of generalized K\"ahler-Einstein structures.
\end{abstract}
\footnote{{\sc Key words}: generalized complex structure, generalized K\"ahler structure, moment map, K\"ahler-Einstein structure, Poisson structure  \\
MSC(2010) 53D18, 53D20, 53D17, 53C26}

\numberwithin{equation}{section}
\section{Introduction}
Let $(X, \ome)$ be a compact symplectic manifold with a symplectic structure $\ome$. 
An almost complex structure $J$ is compatible with $\ome$ if a pair $(J, \ome)$ gives an almost K\"ahler structure on $M$.
We denote by $\til{\mathcal C}_\ome$ the set of almost complex structures which are compatible with $\ome$.
Then $\til{\mathcal C}_\ome$ is an infinite dimensional K\"ahler manifold on which Hamiltonian diffeomorphisms of $(M,\ome)$ act 
$\til{\mathcal C}_\ome$ preserving the K\"ahler structure. Each $J\in \wtil{\mathcal C}_\ome$ gives a Riemannian metric $g(J)$ and
we denote by $s(J)$ the scalar curvature of $g(J)$ which is regarded as a function on $\wtil{\mathcal C}_\ome$.
Then 
the following theorem was established in K\"ahler geometry by Fujiki and Donaldson.
\bgn{theorem}{\text{\rm \cite{Fu_1990}, \cite{Do_1997}}}
The scalar curvature is the moment map on $\til{\mathcal C}_\ome$ for the action of Hamiltonian diffeomorphisms.
\end{theorem}
\noindent The moment map framework in K\"ahler geometry suggests that the existence of constant scalar curvature K\"ahler metrics is inevitably linked with the certain stability in algebraic geometry which leads to well-known Donaldson-Tian-Yau conjecture in K\"ahler geometry. 
\\
\indent Generalized K\"ahler geometry is a successful generalization of ordinary K\"ahler geometry which is equivalent to 
bihermitian geometry satisfying the certain torsion conditions.\\
Many interesting examples of generalized K\"ahler manifolds were already constructed by
holomorphic Poisson structures \cite{Goto_2009}, \cite{Goto_2010}, \cite{Goto_2012}, \cite{Goto_2014}, \cite{Gua_2010}, \cite{Hi_2006},
\cite{Hi_2007}, \cite{Lin_2006}. \\
\indent Main theme of this paper is to pursue an analogue of moment map framework in generalized K\"ahler geometry and to establish the notion of the scalar curvature on a generalized K\"ahler manifold.  
In this paper we assume that a generalized K\"ahler structure consists of commuting two \complex structures $(\J, \J_\psi)$, where $\J$ is an arbitrary almost \complex structure and $\J_\psi$ is induced from a $d$-closed nondegenerate, pure spinor $\psi$ of symplectic type${}^{*1}$.
\footnote{${}^{*1}$ Thus $\psi$ is given by $\psi=e^{b+\sqrt{-1}\ome}$, where $\ome$ is a real symplectic form.
A pair $(\J, \J_\psi)$ is called {\it a generalized K\"ahler structure of symplectic type.}
We can obtain further generalization of moment map framework for any $d$-closed nondegenerate, pure spinor $\psi$. }We construct an invariant function from $\J$ and $\psi$ which is referred to the generalized scalar curvature GR.
Then it turn out that a moment map in generalized K\"ahler geometry is given by the generalized scalar curvature GR.
From the view point of moment map, the notion of generalized Ricci curvature is introduced and the definition of generalized K\"ahler-Einstein structure is provided. 

In order to obtain moment map framework, one may try to follow the same way as in K\"ahler geometry by using the Levi-civita connection and the curvature. However, we need to pave a way without the use of the Levi-civita connection and the curvature in this paper because  
the notion of Levi-civita connection and the curvature in generalized K\"ahler geometry 
are very different from the ones in K\"ahler geometry and are not suitable for our purpose. 
Nondegenerate, pure spinors play a central role rather than \complex structures in this paper.
A nondegenerate, pure spinor is a differential form on a manifold $M$ which induces an almost \complex structure $\J$ by $\ker\phi=E_\J$ and $\ol{\ker\phi}=\ol E_\J$. Conversely, the canonical line bundle $K_\J$ of
a \complex structure $\J$ gives a nondegenerate, pure spinor unique up to multiplication by
nonzero complex functions. 
Let $\{\phi_\a\}$ be trivializations of $K_\J$, where $\phi_\a$ is a nondegenerate, pure spinor which induces \complex structure $\J$. Then the exterior derivative of the differential form $\phi_\a$ is given by 
$$
d\phi_\a=\eta_\a\cdot\phi_\a+ N_\a\cdot\phi_\a
$$
where $\eta_\a$ is a real section $\TT$ and $N_\a$ is also a real section of $\w^3E_\J\oplus \w^3\ol E_\J$. 
A real function $\rho_\a$ is defined by 
$$
\lan\phi_\a,\,\,\, \ol\psi\ran_s=\rho_\a\lan \psi, \,\,\,\ol\psi\ran_s,
$$
where $\lan\,\,,\,\,\ran_s$ denotes the inner metric of Spin representation which is a $2n$-form on $M$.
Then $\J$ acts on $\eta_\a$ by $\J\eta_\a\in \TT$. By Spin representation, $\J\eta_\a$ acts on the differential form $\psi$ by 
$\J\eta_\a\cdot\psi.$ Taking the exterior derivative $d$, we have a differential form 
$d(\J\eta_\a\cdot\psi)$ which is locally defined. We also obtain a differential form 
$d(\J d\log\rho_\a\cdot\psi)$. Then it turns out that 
$-2d(\J\eta_\a\cdot\psi)+d(\J d\log\rho_\a\cdot\psi)$ does not depend on the choice of trivializations of $K_\J$ 
which defines a differential form on $M$ (see Proposition \ref{prop:A differential form}). 
Since $\psi=e^{b+\sqrt{-1}\ome}$, it follows that $d(-2\J\eta_\a+\J d\log\rho_\a)\cdot\ol\psi$ is given by
\bgn{equation}\label{GRic}
 d(-2\J\eta_\a+\J d\log\rho_\a)\cdot\ol\psi=(P-\sqrt{-1}Q)\cdot\ol\psi,
 \end{equation} where $P, Q$ are real $d$-closed $2$-forms.
Thus we define a generalized Ricci form and a generalized scalar curvature GR by 
\bgn{align*}
&\text{\rm GRic}:= -P \quad \text{\rm generalized Ricci form },  \\
&\text{\rm GR} :=n\frac{ P\w\ome^{n-1}}{\ome^n}:\text{\rm generalized scalar curvature}
\end{align*}
where $\ome$ is a symplectic form and GR is a real function (see Definition\ref{def:GRic and GR}).
An almost generalized complex structure $\J$ is compatible with $\psi$ if a pair $(\J, \J_\psi)$ is an almost generalized K\"ahler structure. 
We denote by $\wtil{\A}_\psi(M)$ the set of almost generalized complex structures which are compatible with $\psi$. 
Then it turns out that $\wtil{\A}_\psi(M)$ admits a K\"ahler structure on which the generalized Hamiltonian group 
 acts preserving its K\"ahler structure.
The Lie algebra of generalized Hamiltonian group is given by real smooth functions.

Then our main theorem is the following: \\ \\
\noindent{\sc Theorem} \ref{th:main theorem}\,\,\,\,
There exists a moment map $\mu: {\wtil\A}_\psi(M)\to C^\infty_0(M)^*$ for the generalized Hamiltonian action which is given by the generalized scalar curvature $GR$,
$$
\lan\mu(\J), \, f\ran=(\sqrt{-1})^{-n}\int_M  f(GR_\J) \lan \psi, \,\,\ol\psi\ran_s
$$
\\
In Section 2, we shall give a brief review of almost generalized complex structures focusing on 
nondegenerate, pure spinors and in Section 3, we define an almost generalized K\"ahler structure.  
In Section 4, we recall the stability theorem of generalized K\"ahler structures which is crucial to construct
nontrivial examples of generalized K\"ahler manifolds.
In Section 5, we define a generalized Ricci form GRic and we show that GRic is a representative of the first Chern class of the anticanonical line bundle $K_\J$.
The generalized scalar curvature is obtained from the generalized Ricci form.
The generalized scalar curvature is an invariant function under the action of the extension of volume-preserving diffeomorphisms by
$d$-closed $b$ fields. 
In Section 6, we formulate the moment map framework of generalized K\"ahler geometry. 
After preliminary results are shown in Section 7,
our main theorem is proved in Section 8. In Section 9, we show that infinitesimal deformations of generalized K\"ahler structures with constant generalized scalar curvature are given by an elliptic complex. In particular, the infinitesimal deformations are finite dimensional on a compact manifold.
In Section 10, we give simple
expressions of the generalized Ricci form and the generalized scalar curvature  
of a generalized K\"ahler structure of type $(0,0).$${}^{\dagger}$\footnote{${}^{\dagger}$A generalized K\"ahler structure of type $(0,0)$ corresponds to a degenerate bihermitian structure, i.e, 
$[J_+, J_-]_x\neq0 $ for all $x\in M$.}
 A generalized K\"ahler structure $(\J_\phi, \J_\psi)$ of type $(0,0)$ 
is, by definition, induced from a pair 
$$(\phi=e^{B+\sqrt{-1}\ome_1}, \psi=e^{\sqrt{-1}\ome_2})$$ of $d$-closed, nondegenerate, pure spinors of symplectic types, where $B$ is a real $2$-form and $\ome_1$ and $\ome_2$ are symplectic forms respectively.
Then the $2$-form GRic and the function GR are given by 
\bgn{align*}
\text{GRic}=&-d B\ome_1^{-1}(d\log \frac{\ome_1^n}{\ome_2^n})\\
\text{(GR)}\,\ome_2^n=&\ome_2^{n-1}\w d B\ome_1^{-1}(d\log \frac{\ome_1^n}{\ome_2^n}),
\end{align*}
where $B: T_M\to T^*_M$ and $\ome_i^{-1}: T^*_M\to T_M$ ($i=1,2$).
Then it turns out that the generalized K\"ahler structures coming from hyperK\"ahler structures have vanishing GRic
form.
In Section 11, we define a generalized K\"ahler-Einstein structure. 
In Section 12, we provide nontrivial examples of generalized K\"ahler-Einstein structures which arise as 
Poisson deformations from K\"ahler-Einstein manifolds on which $T^m$ acts preserving its K\"ahler structure.
In particular, the anticanonical divisor of three lines in general position on $\C P^2$ gives a nontrivial example of a
generalized K\"ahler-Einstein structure.

Boulanger obtained remarkable results on the moment map in the cases of  toric generalized K\"ahler manifolds from the view point of 
toric K\"ahler manifolds \cite{Bou_2015}. 
A generalized K\"ahler structure is equivalent to a bihermitian structure with the certain torsion condition. 
From the viewpoint of bihermitian geometry, generalized K\"ahler Ricci flow was introduced \cite{St_2016}. 
Apostolov and Streets discuss Calabi-Yau problem in generalized K\"ahler geometry
\cite{AS_2017}.
It is interesting to find out an expression of our moment map in terms of bihermitian geometry.
There is a remarkable link between generalized geometry and noncommutative algebraic geometry. 
It is quite natural to ask whether the existence of generalized K\"ahler structure with constant generalized scalar curvature is related with a stability on a  noncommutative algebraic manifold.

{\bf Acknowledgement} The author would like to thank Professor N. J. Hitchin for his interests and 
remarkable comments on his results. 
The author also thanks Shinnosuke Okawa for valuable discussions on the relation between generalized geometry and
noncommutative algebraic geometry.
The author would like to thank Marco Gualtieri and Vesti Apostolov for valuable comments.
\section{Generalized complex structures}
Let $M$ be a differentiable manifold of real dimension $2n$.
The bilinear form $\lan\,\,,\,\,\ran_{\scriptstyle T\oplus T^*}$ on 
the direct sum $T_M \oplus T^*_M$ over a differentiable manifold $M$ of dim$=2n$ is defined by 
$$\lan v+\xi, u+\eta \ran_{\tt}=\frac12\(\xi(u)+\eta(v)\),\quad  u, v\in T_M, \xi, \eta\in T^*_M .$$
Let SO$(\TT)$ be the fibre bundle over $M$ with fibre SO$(2n, 2n)$ which is 
a subbundle of End$(\TT)$  preserving the bilinear form $\lan\,,\,\ran_s$ 
 An almost \complex structure $\J$ is a section of SO$(\TT)$ satisfying $\J^2=-\id.$ Then as in the case of almost complex structures, an almost \complex structure $\J$ yields the eigenspace decomposition :
$(T_M \oplus T^*_M)^\C =E_\J \oplus \ol E_\J$, where 
$E_\J$ is $-\sqrt{-1}$-eigenspaces and  $\ol{E}_\J$ is the complex conjugate of $E_\J$. 
The Courant bracket of $\TT$ is defined by 
$$
 [u+\xi, v+\eta]_{\text{cou}}=[u,v]+{\mathcal L}_u\eta-{\mathcal L}_v\xi-\frac12(di_u\eta-di_v\xi),
 $$
 where $u, v\in TM$ and $\xi, \eta$ is $T^*M$.
If $E_\J$ is involutive with respect to the Courant bracket, then $\J$ is a generalized complex structure, that is, $[e_1, e_2]_{\text{cou}}\in \Gam(E_\J)$  for any two elements 
 $e_1=u+\xi,\,\, e_2=v+\eta\in \Gam(E_\J)$.

Let $\CL(T_M \oplus T^*_M)$ be the Clifford algebra bundle which is 
a fibre bundle with fibre the Clifford algebra $\CL(2n, 2n)$ with respect to $\lan\,,\,\ran_{\scriptstyle T \oplus T^*}$ on $M$.
Then a vector $v$ acts on the space of differential forms $\oplus_{p=0}^{2n}\w^pT^*M$ by 
the interior product $i_v$ and a $1$-form acts on $\oplus_{p=0}^{2n}\w^pT^*M$ by the exterior product $\theta\wedge$, respectively.
Then the space of differential forms gives a representation of the Clifford algebra $\CL(\TT)$ which is 
the spin representation of $\CL(\TT)$. 
Thus
the spin representation of the Clifford algebra arises as the space of differential forms $$\w^\bullet T^*_M=\oplus_p\w^pT^*_M=\w^{\even}T^*_M\oplus\w^{\odd}T^*_M.$$ 
The inner product $\lan\,,\,\ran_s$ of the spin representation is given by 
$$
\lan \a, \,\,\,\beta\ran_s:=(\a\w\sig\beta)_{[2n]},
$$
where $(\a\w\sig\beta)_{[2n]}$ is the component of degree $2n$ of $\a\w\sig\beta\in\oplus_p \w^pT^*M$ and 
$\sig$ denotes the Clifford involution which is given by 
$$
\sig\beta =\bgn{cases}&+\beta\qquad \deg\beta \equiv 0, 1\,\,\mod 4 \\ 
&-\beta\qquad \deg\beta\equiv 2,3\,\,\mod 4\end{cases}
$$

We define $\ker\Phi:=\{ e\in (T_M\oplus T^*_M)^\C\, |\, e\cdot\Phi=0\, \}$ for a differential form $\Phi
\in \w^{\even/\odd}T^*_M.$
If $\ker\Phi$ is maximal isotropic, i.e., $\dim_\C\ker\Phi=2n$, then $\Phi$ is called {\it a pure spinor} of even/odd type.

A pure spinor $\Phi$ is {\it nondegenerate} if $\ker\Phi\cap\ol{\ker\Phi}=\{0\}$, i.e., 
$(T_M\oplus T^*_M)^\C=\ker\Phi\oplus\ol{\ker\Phi}$.
Then a nondegenerate, pure spinor $\Phi\in \w^\bullet T^*_M$ gives an almost generalized complex structure $\J_{\Phi}$ which satisfies 
$$
\J_\Phi e =
\bgn{cases}
&-\sqrt{-1}e, \quad e\in \ker\Phi\\
&+\sqrt{-1}e, \quad e\in \ol{\ker\Phi}
\end{cases}
$$
Conversely, an almost \complex structure $\J$ locally arises as $\J_\Phi$ for a nondegenerate, pure spinor $\Phi$ which is unique up to multiplication by
non-zero functions.  Thus an almost \complex structure yields the canonical line bundle $K_{\J}:=\C\lan \Phi\ran$ which is a complex line bundle locally generated by a nondegenerate, pure spinor $\Phi$ satisfying 
$\J=\J_\Phi$.
An \complex structure 
$\J_\Phi$ is integrable if and only if $d\Phi=\eta\cdot\Phi$ for a section $\eta\in T_M\oplus T^*_M$. 
The {\it type number} of $\J=\J_\Phi$ is defined as the minimal degree of the differential form $\Phi$. Note that type number Type $\J$ is a function on a manifold which is not a constant in general.
\bgn{example}
Let $J$ be a complex structure on a manifold $M$ and $\J^*$ the complex structure on the dual  bundle $T^*M$ which is given by $J^*\xi(v)=\xi (Jv)$ for $v\in TM$ and $\xi\in T^*M$.
Then a \complex structure $\J_J$ is given by the following matrix
$$\J_J=\bgn{pmatrix}J&0\\0&-J^*
\end{pmatrix},$$
Then the canonical line bundle is the ordinary one which is generated by complex forms of type $(n,0)$.
Thus we have  Type $\J_J =n.$
\end{example}
\bgn{example}
Let $\ome$ be a symplectic structure on $M$ and $\h\ome$ the isomorphism from $TM$ to $T^*M$ given by $\h\ome(v):=i_v\ome$. We denote by $\h\ome^{-1}$ the inverse map from $T^*M$ to $TM$.
Then a \complex structure $\J_\psi$ is given by the following
$$\J_\psi=\bgn{pmatrix}0&-\h\ome^{-1}\\
\h\ome&0
\end{pmatrix},\quad\text{\rm Type $\J_\psi =0$}$$
Then the canonical line bundle is given by  the differential form $\psi=e^{\sqrt{-1}\ome}$. 
Thus Type $\J_\psi=0.$
\end{example}
\bgn{example}[$b$-field action]
A $d$-closed $2$-form $b$ acts on a \complex structure by the adjoint action of Spin group $e^b$ which provides
a \complex structure $\Ad_{e^b}\J=e^b\circ \J\circ e^{-b}$. 
\end{example}
\bgn{example}[Poisson deformations]
Let $\beta$ be a holomorphic Poisson structure on a complex manifold. Then the adjoint action of Spin group $e^\beta$ gives deformations of new \complex structures by 
$\J_{\beta t}:=\Ad_{\beta^{Re} t}\J_J$.  Then Type ${\J_{\beta t}}_x=n-2$ rank of $\beta_x$ at $x\in M$,
which is called the Jumping phenomena of type number.
\end{example}
Let $(M, \J)$ be a generalized complex manifold and $\ol E_\J$ the eigenspace of eigenvalue $\sqrt{-1}$.
Then we have the Lie algebroid complex $\w^\bullet\ol{E}_\J$:
$$
0\arrow\w^0\ol E_\J\overset{\ol\pa_\J}\arrow\w^1\ol E_\J\overset{\ol\pa_\J}\arrow\w^2\ol E_\J\overset{\ol\pa_\J}\arrow\w^3\ol E_\J\arrow\cdots
$$
The Lie algebroid complex is the deformation complex of \complex structures. 
In fact, $\e\in \w^2\ol E_\J$ gives deformed isotropic subbundle 
$E_\e:=\{ e+[\e, e]\, |\, e\in E_\J\}$. 
Then $E_\e$ yields deformations of \complex structures if and only if $\e$ satisfies Generalized Mauer-Cartan equation
$$
\ol{\pa}_\J\e+\frac12[\e, \e]_{\Sch}=0,
$$
where $[\e, \e]_{\Sch}$ denotes the Schouten bracket. 
The Kuranishi space of generalized complex structures is constructed.

Then the second cohomology group $H^2(\w^\bullet\ol E_\J)$ of the Lie algebraic complex gives the infinitesimal deformations of \complex structures and the third one 
$H^3(\w^\bullet\ol E_\J)$ is the obstruction space to deformations of \complex structures.

Let $\{e_i\}_{i=1}^n$ be a local basis of $E_\J$ for an almost \complex structure $\J$, 
where $\lan e_i, \ol e_j\ran_{\tt}=\del_{i,j}$.
The the almost \complex structure $\J$ is written as an element of Clifford algebra,
$$
\J=\frac{\sqrt{-1}}2\sum_i e_i\cdot\ol {e}_i,
$$
where $\J$ acts on $\TT$ by the adjoint action $[\J, \,]$. 
Thus we have $[\J, e_i]=-\sqrt{-1}e_i$ and $[\J, \ol e_i]=\sqrt{-1}e_i$.
An almost \complex structure $\J$ acts on differential forms by the Spin representation which gives the decomposition:
\bgn{equation}
\w^\bullet T^*_M=U^{-n}\oplus U^{-n+1}\oplus\cdots U^{n}
\end{equation}
\section{Almost generalized K\"ahler structures}

\bgn{definition}

{\it An almost generalized K\"ahler structure} is a pair $(\J_1, \J_2)$ consisting of two commuting almost \complex structures 
$\J_1, \J_2$ such that $\h G:=-\J_1\circ\J_2=-\J_2\circ \J_1$ gives a positive definite symmetric form 
$G:=\lan \h G\,\,,  \,\,\ran$ on $T_M\oplus T_M^*$, 
We call $G$ {\it a generalized metric}.
{\it A generalized K\"ahler structure } is an almost \GK structure $(\J_1, \J_2)$ such that both $\J_1$ and $\J_2$ are \complex structures. 
\end{definition}
$\J_i$ gives the decomposition $(\TT)^\C=E_{\J_i}\oplus\ol E_{\J_i}$ for $i=1,2$.
Since $\J_1$ and $\J_2$ are commutative, we have the simultaneous eigenspace decomposition 
$$
(\TT)^\C=(E_{\J_1}\cap E_{\J_2})\oplus (\ol E_{\J_1}\cap \ol E_{\J_2})\oplus (E_{\J_1}\cap \ol E_{\J_2})\oplus
(\ol E_{\J_1}\cap E_{\J_2}).
$$
Since $\h G^2=+\id$,
The generalized metric $\h G$ also gives the eigenspace decomposition: $\TT=C_+\oplus C_-$, 
where $C_\pm$ denote the eigenspaces of $\h G$ of eigenvalues $\pm1$. 
We denote by $E_{\J_1}^\pm$ the intersection $E_{\J_1}\cap C^\C_\pm$. 
Then it follows 
\bgn{align*}
&E_{\J_1}\cap E_{\J_2}=E_{\J_1}^+,  \quad \ol E_{\J_1}\cap \ol E_{\J_2}=\ol E_{\J_1}^+\\
&E_{\J_1}\cap \ol E_{\J_2}=E_{\J_1}^-,\quad \ol E_{\J_1}\cap E_{\J_2}=\ol E_{\J_1}^-
\end{align*}

\bgn{example}
Let $X=(M, J,\ome)$ be a K\"ahler manifold. Then the pair $(\J_J, \J_\psi)$ is a generalized K\"ahler where 
$\psi=\exp(\sqrt{-1}\ome)$.
\end{example}
\section{The stability theorem of generalized K\"ahler manifolds}
It is known that the stability theorem of ordinary K\"ahler manifolds holds
\bgn{theorem}[Kodaira-Spencer]
Let $X=(M,J)$ be a compact K\"ahler manifold and $X_t$ small deformations of $X=X_0$ as complex manifolds.
Then $X_t$ inherits a K\"ahler structure. 
\end{theorem}

The following stability theorem of generalized K\"ahler structures shows that 
there are many intriguing examples
of generalized K\"ahler manifolds of symplectic type.
\bgn{theorem}{\rm \cite{Goto_2010}}
Let $X=(M,J,\ome)$ be a compact K\"ahler manifold and $(\J_J, \J_\psi)$ the induced generalized K\"ahler structure, 
where $\psi=e^{\sqrt{-1}\ome}$. 
If there are analytic deformations $\{\J_t\}$ of $\J_0=\J_J$ as \complex structures, then there are deformations of $d$-closed nondegenerate, pure spinors $\{\psi_t\}$ such that 
pairs $(\J_t, \J_{\psi_t})$ are generalized K\"ahler structures, where $\psi_0=\psi$
\end{theorem}
\section{Generalized Ricci curvature and generalized scalar curvature}
We use the same notation as before.
Let $\J$ be an almost complex structure on $M$ with 
trivializations $\{\phi_\a\}$ of the canonical line bundle $K_\J.$
Then recall that $\eta_\a$ is given by 
\bgn{equation}\label{eta}
d\phi_\a=\eta_\a\cdot\phi_\a+ N_a\cdot\phi_\a,  
\end{equation}
where $\eta_\a\in T_M\oplus T^*_M$ and $N_\a\in \w^3E_\J\oplus\w^3\ol E_\J$ are real sections, i.e., 
$\eta_a=\ol\eta_\a, \,\,N_\a=\ol N_\a$. 
Because of the reality condition, $\eta_\a$ and $N_\a$ are uniquely determined. 
Let $(\J, \psi)$ be an almost generalized K\"ahler structure of symplectic type.
Then recall that a real function $\rho_\a$ on $U_\a$ is given by 
\bgn{equation}\label{rho}
\lan \phi_a, \,\,\ol\phi_\a\ran_s=\rho_\a\lan \psi, \ol\psi\ran_s
\end{equation}

\bgn{proposition}\label{prop:A differential form}
A differential form 
$d(-2\J\eta_\a+\J d\log{\rho_\a})\cdot\ol\psi$ does not depend on the choice of trivializations $\{\phi_\a\}$ of $K_\J$.
\end{proposition}
\bgn{proof}
Let $e^{\kappa_{\a,\beta}}$ be the transition function on the intersection $U_\a\cap U_\beta$. 
Then we have 
$\phi_\a=e^{\kappa_{\a,\beta}}\phi_\beta$.
Since $d\phi_\beta=(\eta_\beta+N_\beta)\cdot\phi_\beta$, we have
\bgn{align*}
d\phi_\a=&d(e^{\kappa_{\a,\beta}}\phi_\beta)=
d\kappa_{\a,\beta}\cdot e^{\kappa_{\a,\beta}}\phi_\beta+e^{\kappa_{\a,\beta}}(\eta_\beta+N_\beta)\phi_\beta\\
=&d\kappa_{\a,\beta}\cdot \phi_\a+(\eta_\beta+N_\beta)\cdot\phi_\a
\end{align*}
Thus we have $(\eta_\a+N_\a)\cdot\phi_\a=(\eta_\beta+d\kappa_{\a,\beta}+N_\beta)\cdot\phi_\a $.
Since $N_\a, N_\beta\in \w^3 E_\J\oplus \w^3\ol E_\J$ and $\eta_a, \eta_b ,\ d\kappa_{\a,\beta}\in \TT$, we have $N_\a=N_\beta$ and $\eta_\a\cdot\phi_\a=(\eta_\beta+d\kappa_{\a,\beta})\cdot\phi_a$.
Since $\eta_\a, \eta_\beta$ are real, it follows that we have 
\bgn{equation}\label{eta henkan}
\eta_\a=\eta_\beta+\ol\pa_\J\kappa_{\a,\beta}+\pa_\J  \ol{\kappa_{\a,\beta}}
\end{equation}
We also have 
\bgn{equation}\label{rho henkan}
\rho_\a=\rho_\beta e^{\kappa_{\a,\beta}+\ol\kappa_{\a,\beta}}
\end{equation}
Since $d((d\kappa_{\a,\beta})\cdot\psi)=0$, it follows from $d\kappa_{\a,\beta}=\pa_\J \kappa_{\a,\beta}+\ol\pa_\J\kappa_{\a,\beta}$ that we have
\bgn{equation}\label{key 1}
d(\pa_\J \kappa_{\a,\beta})\cdot\psi+d(\ol\pa_J\kappa_{\a,\beta})\cdot\psi=0
\end{equation}
Since $d((d\ol\kappa_{\a,\beta})\cdot\psi)=0$, we also have 
\bgn{equation}\label{key 2}
d(\ol\pa_\J \ol\kappa_{\a,\beta})\cdot\psi+d(\pa_J\ol\kappa_{\a,\beta})\cdot\psi=0
\end{equation}
Applying (\ref{eta henkan}), (\ref{key 1}) and (\ref{key 2}),
 we have 
\bgn{align*}
&d(-2\J\eta_\a+\J d\log\rho_\a)\cdot\psi-d(-2\J\eta_\beta+\J d\log\rho_\beta)\cdot\psi\\
=&-2d\J(\ol\pa_\J\kappa_{\a,\beta}+\pa_\J  \ol{\kappa_{\a,\beta}})\cdot\psi+d\J d(\kappa_{\a,\beta}+\ol\kappa_{\a,\beta})\cdot\psi\\
=&-2\sqrt{-1}d(\ol\pa_\J\kappa_{\a,\beta}-\pa_\J  \ol{\kappa_{\a,\beta}})\cdot\psi+\sqrt{-1}d(\ol\pa_\J \kappa_{\a,\beta}-
\pa_\J\kappa_{\a,\beta}+\ol\pa_\J \ol\kappa_{\a,\beta}-\pa_\J\ol\kappa_{\a,\beta})\cdot\psi\\
=&-2\sqrt{-1}d\ol\pa_\J(\kappa_{\a,\beta}+ \ol{\kappa_{\a,\beta}})\cdot\psi+2\sqrt{-1}d\ol\pa_\J(\kappa_{\a,\beta}+ \ol{\kappa_{\a,\beta}})\cdot\psi=0
\end{align*}
Thus we have the result${}^{*2}$\footnote{${}^{*2}$In this proof, note that we do not use the integrability of $\J$.}.
\end{proof}
Hence $d(-2\J\eta_\a+\J d\log{\rho_\a})\cdot\ol\psi$ yields a globally defined differential form on $M$.
Since $\psi=e^{b+\sqrt{-1}\ome}$, it follows that $d(-2\J\eta_\a+\J d\log\rho_\a)\cdot\ol\psi$ is given by
\bgn{equation}\label{GRic}
 d(-2\J\eta_\a+\J d\log\rho_\a)\cdot\ol\psi=(P-\sqrt{-1}Q)\cdot\ol\psi,
 \end{equation} where $P, Q$ are real $d$-closed $2$-forms.
 In fact, $-2\J\eta_\a+\J d\log{\rho_\a}$ is written as $v+\theta\in \TT$ for a vector $v$ and a $1$-form $\theta$ and 
 then $-2\J\eta_\a+\J d\log{\rho_\a}\cdot\ol\psi$ is given by 
 $(i_vb-\sqrt{-1}i_v\ome+\theta)\w\ol\psi$. Thus $P$ and $Q$ are given by 
 $P=di_v b+d\theta$ and $Q=di_v\ome$.
\bgn{remark}
Since $N_\a=N_\beta$, we have a globally defined section $N\in \w^3 E_\J\oplus \w^3\ol E_\J$ which is 
the Nijenhuis type tensor, that is, $\J$ is integrable if and only if $N=0$.
\end{remark} 
\bgn{definition}\label{def:GRic and GR}{[Generalized Ricci form and generalized scalar curvature]}
We define a generalized Ricci form GRic to be a $d$-closed $2$-form $P$ in (\ref{GRic}) and  
we define a generalized scalar curvature GR to be a real function on $M$ which is given by the following,
\bgn{align*}
&\text{\rm GRic}:=- P \quad \text{\rm generalized Ricci form },  \\
&\text{\rm GR} :=\frac{n P\w\ome^{n-1}}{\ome^n}:\text{\rm generalized scalar curvature}
\end{align*}
where $\ome$ is a symplectic form.
\end{definition}

A diffeomorphism $F$ of $M$ acts on $(\J, \psi)$ to give an almost generalized K\"ahler structure 
$(\J', \psi')$.
We denote by GR' generalized scalar curvature of $(\J', \psi').$
Then we have 
\bgn{proposition}
 $$
\GR'=F^*(\GR),
$$
that is , GR is equivalent under the action of diffeomorphisms.
Further GR is invariant under the action of $d$-closed $b-$fields.
\end{proposition}
\bgn{proof}
A diffeomorphism $F$ of $M$ induces the bundle map $F_\#$ of $\TT$ by 
$F_\#(v+\theta)=F^{-1}_*(v)+ F^*\theta$ for $v\in T_M$ and $\theta\in T\*_M.$
Then we see that 
$F_\#(v+\theta)\cdot F^*(\a)=F^*((v+\theta)\cdot\a)$ for a differential form $\a.$
Let $b$ be a real $d$-closed $2$-form.
Then $e^b$ is regarded as an element of Spin group of the Clifford algebra of $\TT$ which acts on 
 differential forms by the wedge produce of $e^b $. 
Then we have the adjoint action $\Ad_{e^b}$ on $\TT$ by 
$\Ad_{e^b}(v+\theta) :=e^b (v+\theta) e^{-b}=v- i_vb + \theta.$
Then we see that 
$$(\J', \psi')=(F_\#\circ\J \circ F_\#^{-1}, \,\, F^*\psi).$$
Then it follows that 
$\phi_\a'=F^*\phi_\a$ is the nondegenerate pure spinor which induces $\J'.$
We define $\eta_\a'$ by $d\phi_\a'=\eta_\a'\cdot\phi_\a'$.
Thus we  have
$$
d\phi_\a'= F^*d\phi_\a = F^*(\eta_\a\cdot \phi_\a)=F_\#(\eta_\a)\cdot\phi_\a'
$$
Thus we see that $\eta_\a'=F_\#(\eta_\a).$ 
The function $\rho'_\a$ is given by 
$$
\lan \phi_\a', \,\,\ol\phi_\a'\ran_s=\rho_\a'\lan \psi', \,\,\ol\psi'\ran_s
$$
 Thus we have $\rho_\a'=F^*\rho_\a.$
 Then we see 
 \bgn{align}
 (\J'\eta_\a')\cdot\psi'=& F_\#\circ \J \circ F_\#^{-1}(F_\#(\eta_\a))\cdot F^*\psi\\
 =&F_\#(\J\eta_\a)\circ F^*\psi \\
 =&F^*(\J\eta_\a\cdot\psi)
 \end{align} 
 We also have 
 \bgn{align}
 \J' (d\log\rho_\a')\cdot\psi'=&F_\#\circ\J\circ F_\#^{-1}(F^*(d\log\rho_\a)F^*{\psi)}\\
 =&F_\#(\J(d\log\rho_\a))\cdot F^*\psi\\
 =&F^*(\J(d\log\rho_\a)\cdot\psi)
 \end{align}
 Thus we obtain 
 \bgn{equation}\label{eq:d(-2J'}
 d(-2\J'\eta_\a'+\J' d\log\rho_\a')\cdot\ol\psi'=F^*\(d(-2\J\eta_\a+\J d\log\rho_\a)\cdot\ol\psi\)
 \end{equation}
 Since GR is given by the real part of the following:
 $$GR:=\Re \frac{\sqrt{-1}}2\frac{\lan \psi, \,\,\, d(-2\J\eta_\a+\J d\log{\rho_\a})\cdot\ol\psi\ran_s}
{\lan \psi, \,\,\,\ol\psi\ran_s}$$
From (\ref{eq:d(-2J'}), we have GR'$=F^*(\GR).$

We denote by $(\J_b, \psi_b)$ the pair given by the action of $e^b$ on $(\J, \psi)$.
Then $\J_b$ is induced from $e^b\cdot \phi_\a$ and $\psi_b=e^b\cdot\psi.$
Thus we have $\eta_\a^b =\Ad_{e^b} (\eta_\a)$ and $\rho_\a^b =\rho_\a.$ 
Then we see that 
$$
 d(-2\J_b\eta_\a^b+\J_b d\log\rho_\a^b)\cdot\ol\psi_b=e^b\(d(-2\J\eta_\a+\J d\log\rho_\a)\cdot\ol\psi\)
$$
Since $\lan\,\,,\,\,\ran_s$ is invariant under the action of $e^b,$
we see that GR is invariant under the action of $e^b.$
\end{proof}
We denote by [GRic] the cohomology class of a real $d$-closed $2$-form GRic. Then we have
\bgn{proposition}  The cohomology class  $[\text{\rm GRic}]$ is given by the $1$-st Chern class, 
$$
[\text{\rm GRic}]=4\pi c_1(K_\J^{-1})\in H^2(M)
$$
\end{proposition}
\bgn{proof}
We calculate the spectral sequence from de Rham to \v Ceck cohomology to determine
a representative of {\v C}eck cohomology group given by $d$-closed form GRic. 
 $d(-2\J\eta_\a+\J d\log\rho_\a)\cdot\psi$ is $d$-exact on $U_\a$.
On $U_\a\cap U_\beta$, it follows from (\ref{eta henkan}) and (\ref{rho henkan}) that we have
\bgn{align*}
&(-2\J\eta_\a+\J d\log\rho_\a)\cdot\psi-(-2\J\eta_\beta+\J d\log\rho_\beta)\cdot\psi\\
=&(-2\J(\eta_\a-\eta_\beta)+\J d (\kappa_{\a,\beta}+\ol \kappa_{\a,\beta}))\cdot\psi\\
=&-2\J(\ol \pa_\J\kappa_{\a,\beta}+\pa \ol \kappa_{\a,\beta})+\J (\ol\pa_\J \kappa_{\a,\beta}+\pa_\J \kappa_{\a,\beta}+\ol\pa_\J\ol \kappa_{\a,\beta}+\pa_\J \ol\kappa_{\a,\beta})\cdot\psi\\
=&-\sqrt{-1} d(\kappa_{\a,\beta}-\ol\kappa_{\a\beta})\cdot\psi\\
=&2 d k_{\a,\beta}^{\Im}\cdot\psi,
\end{align*}
where $\kappa_{\a,\beta}^{\Im}$ denotes the imaginary part of $\kappa_{\a,\beta}$. 
Thus we have a \v Ceck representative, 
$$
2 (\kappa_{\a,\beta}^{\Im}+\kappa_{\beta,\gam}^{\Im}+\kappa_{\gam, \a}^{\Im})\cdot\psi
$$
Thus the representative of the class $[P]$ is given by 
$2 (\kappa_{\a,\beta}^{\Im}+\kappa_{\beta,\gam}^{\Im}+\kappa_{\gam, \a}^{\Im}).$

The $1$-st Chern class $c_1(K_\J)$ has a \v Ceck representative 
$$
c_{\a,\beta,\gam}=\frac 1{2\pi\sqrt{-1}}(\kappa_{\a,\beta}+\kappa_{\beta,\gam}+\kappa_{\gam, \a})
\equiv\frac1{2\pi}(\kappa_{\a,\beta}^{\Im}+\kappa_{\beta,\gam}^{\Im}+\kappa_{\gam, \a}^{\Im})
$$
Thus we have $[P] =4\pi c_1(K_\J)$.
Since GRic $=-P$, we  obtain the result.

\end{proof}

\bgn{example}
A GK structure $(\J_J, \psi=e^{\sqrt{-1}\ome})$ is induced from the genuine K\"ahler structure. 
Then GRic and GR are the ordinary Ricci curvature and scalar curvature, respectively. 
In fact, we have $\phi_\a$ to be a holomorphic $n$ form 
$\phi_\a=dz_1\w\cdots\w dz_n$ and $\psi=e^{\sqrt{-1}\ome}$ and
$\lan\Ome_\a,\ol \Ome_\a\ran_S=\rho_\a\lan\psi,\ol\psi\ran_S$.
Thus $d\J d\log\rho_\a=-2\sqrt{-1}\pa\ol\pa\log\det g_{i,\ol{j}}$  is the ordinary Ricci form. 
\end{example}
\bgn{remark}
We can generalize our construction of  GR to the cases where $\psi$ is an arbitrary $d$-closed, nondegenerate, pure spinor. 
In fact, $ d(-2\J\eta_\a+\J d\log{\rho_\a})\cdot\ol\psi$ is still a representative of the first Chern class of $K_\J$ together with the class $[\ol\psi]$ and 
$$GR^\C:= \frac{\sqrt{-1}}2\frac{\lan \psi, \,\,\,, d(-2\J\eta_\a+\J d\log{\rho_\a})\cdot\ol\psi\ran_s}
{\lan \psi, \,\,\,\ol\psi\ran_s}$$
is an equivalent complex function under the action of diffeomorphisms which is invariant under the action of $d$-closed $b$-fields. In this general case, we define GR to be the real part of GR$^\C$. 
Then we have 
\bgn{align*}
(\sqrt{-1})^{-n}(GR)\lan \psi, \,\,\,\ol\psi\ran_s=\Re(\sqrt{-1})^{-n+1}
\lan \psi, \,\,\,d(-\J \eta_\a+\frac12 \J d\log\rho_\a)\cdot\ol\psi \ran_s,
\end{align*}
where Re stands for the real part. 
The real part is also written as 
\bgn{align}\label{GR in general} 
(\sqrt{-1})^{-n}(GR)\lan \psi, \,\,\,\ol\psi\ran_s=&c_n \lan \psi, \,\,\,d(-\J \eta_\a+\frac12 \J d\log\rho_\a)\cdot\ol\psi \ran_s\\
-&c_n\lan d(-\J\eta_\a+\frac12\J d\log\rho_\a)\cdot\psi, \,\,\,\ol\psi\ran_s,\notag
\end{align}
where $c_n=\frac12(\sqrt{-1})^{-n+1}$.
\end{remark}
\bgn{example}[generalized Calabi-Yau metrical structure]
If a generalized K\"ahler structure is induced from a pair $(\phi, \psi)$ which consists of $d$-closed, nondegenerate, pure spinors
such that $\lan\phi,\ol\phi\ran_S=\lan \psi, \ol\psi\ran_S$, then it is called {\it a generalized Calabi-Yau metrical structure.} 
Since $\rho_\a=1$ and $\eta_\a=0$, 
it follows that we have GR$^\C=0$.
\end{example}

\section{Generalized scalar curvature as moment map}
Let ${{\mathcal GC}}(M)$ be the set of \complex structures on a differentiable compact manifold $M$ of dimension $2n$,
 that is, 
$${{\mathcal GC}}(M):=\{ \J\,:\,\text{\rm \complex structure on }M\,\}.$$
We denote by ${{\mathcal GK}}(M)$ the set of generalized K\"ahler structures on $M$, that is,
$${{\mathcal GK}}(M):=\{(\J_0, \J_1):\,\text{\rm generalized K\"ahler structure on }M\, \}.$$
We also define $\wtil{{\mathcal GC}}(M)$ as the set of almost \complex structures on $M$,
$$\wtil{{\mathcal GC}}(M):=\{ \J\, :\text{\rm almost \complex structure on }M\,\}.$$
We denote by $\wtil{{\mathcal GK}}(M)$ the set of almost generalized K\"ahler structures, 
$$\wtil{{\mathcal GK}}(M):=\{(\J_0, \J_1):\,\text{\rm almost generalized K\"ahler structure on }M\, \}.$$
Let $\psi$ be a $d$-closed, non-degenerate, pure spinor which induces $\J_\psi$.
The spinor inner product of $\psi$ is given by
$\lan\psi, \ol\psi\ran_S=(\phi\w\sig\psi)_{[2n]}.$ 
In particular, if $\psi:=e^{b+\frac{\sqrt{-1}}{2}\ome}$, then we have the volume form
$$\lan\psi, \ol\psi\ran_S=\frac{(\sqrt{-1})^n}{n!}\ome^n.$$
An almost \complex structure $\J$ is {\it $\psi$-compatible } if and only if the pair $(\J, \J_\psi)$
is an almost generalized K\"ahler structure.
Let ${{\A}}_\psi(M)$ be the set of $\psi$-compatible \complex structure, that is  
$${{\A}}_\psi(M):=\{\, \J\in {{\mathcal GC}}\, :\, (\J,\J_\psi)\in {{\mathcal GK}}\, \}.$$
We also define $\wtil{{\A}}_\psi(M)$ to be the set of $\psi$-compatible almost \complex structures,
 $$\wtil{{\A}}_\psi(M):=\{\, \J\in {\wtil{\mathcal GC}}\, :\, (\J,\J_\psi)\in {\wtil{\mathcal GK}}\, \}.$$

For each point $x\in M$, we define $\wtil{\A}_\psi(M)_x$ to be the set of $\psi_x$-compatible almost \complex structures , that is, 
 $${\wtil{\A}}_\psi(M)_x:=\{\, \J_x\, |(\J_x, \J_{\psi, x}): \text{\rm almost generalized K\"ahler structure at } x \, \}.$$
 Then we see that $\wtil{\A}_\psi(M)_x$ is given by the Riemannian Symmetric space of type  AIII${}^{\dag}$
$$U(n,n)/U(n)\times U(n)$$ which is biholomorphic to the complex bounded domain 
 $\{\, h\in M_n(\C)\, |\, 1_n-h^*h>0\, \},$ where $M_n(\C)$ denotes the set of complex matrices of $n\times n.$ 
 \footnote{${}^{\dag}$ In K\"ahler geometry, the set of almost complex structures compatible with a symplectic structure $\ome$ is given by 
 the Riemannian symmetric space Sp$(2n)/U(n)$ which is biholomorphic to the Siegel upper half plane 
 $\{\, h\in \GL_n(\C)\, |\, 1_n-h^*h>0,\, h^t=h\}$. \,}
 Let $P_\psi$ be the fibre bundle over $M$ with fibre ${\wtil{\A}}_\psi(M)_x$, that is, 
$$P_\psi:=\bigcup_{x\in M}{\wtil{\A}_\psi(M)_x}\to M,$$
Then $\wtil{\A}_\psi(M)$ is given by sections $\Gam (M, P_\psi)$ which contains  ${{\A}}_\psi(M)$. 
We can introduce a Sobolev norm on $\wtil{\A}_\psi(M)$ such that $\wtil{\A}_\psi(M)$ becomes a Banach manifold in the usual way.
The tangent bundle of $\wtil{\A}_\psi(M)$ at $\J$ is given by 
$$T_{\J}\wtil{\A}_\psi(M)=\{\, \dot{\J}\in\text{\rm so}(T_M\oplus T^*_M)\,:\, \dot{\J}\J+\J\dot{\J}=0,\, \dot{\J}\J_\psi=\J_\psi\dot{\J}\, \},$$
where so$(\TT)$ denotes the set of sections of Lie algebra bundle of SO$(\TT)$.
Then it follows that there exists an almost complex structure $J_{\wtil{\A}_\psi}$ on $\wtil{\A}_\psi(M)$which is given by 
$$
J_{\til{\A}_\psi}(\dot{\J}):=\J\dot{\J}, \qquad (\,\,\dot{\J}\in T_{\J}\wtil{\A}_\psi(M) \,\,)
$$
We also have a Riemannian metric $g_{\wtil{\A}_\psi}$ and a $2$-form $\ome_{\wtil{\A}_\psi}$ on ${\wtil\A}_\psi(M)$ by 
\bgn{align}\label{Apsi}
&g_{\wtil{\A}_\psi}(\dot{\J_1},\dot{\J_2}):=\frac1{(\sqrt{-1})^n}\int_M \tr(\dot{\J_1}\dot{\J_2})\lan\psi,\,\ol\psi\ran_S
\\
&\ome_{\wtil{\A}_\psi}(\dot{\J_1},\dot{\J_2}):=\frac{-1}{(\sqrt{-1})^n}\int_M \tr(\J\dot{\J_1}\dot{\J_2})\lan\psi,\,\ol\psi\ran_S
\end{align}
for $\dot{\J_1}, \dot{\J_2}\in T_{\J}\wtil{\A}_\psi(M)$. 
\bgn{proposition}
$J_{\til{\A}_\psi}$ is an integrable almost complex structure on $\wtil{\A}_\psi(M)$ and
$\ome_{\wtil{\A}_\psi}$ is a K\"ahler form on $\wtil{\A}_\psi(M).$
\end{proposition}
\bgn{proof}
Let $\J_{V}$ be an almost generalized complex structure on a real vector space $V$ of dimension $2n$.
We denote by $X_n$ the Riemannian symmetric space  
$U(n,n)/U(n)\times U(n)$ which is identified with the set of almost generalized complex structures
compatible with $\J_{V}.$
We already see that $\wtil{\A}_\psi(M)$ is the set of global sections of the fibre bundle $P_\psi$ over a manifold $M$ with fibre 
$X_n$ which is biholomorphic to the bounded domain $\{\, h\in M_n(\C)\, |\, 1_n-h^*h>0\, \}.$
If $\wtil{\A}_\psi(M)$ is not empty, we have a global section $\J_0$.
Then the fibre bundle is identified with the space of maps 
from $M$ to the complex bounded domain $\{\, h\in M_n(\C)\, |\, 1_n-h^*h>0\, \}$ which is open set in the complex vector space $M_n(\C)$. 
Since the almost complex structure $J_{\til{\A}_\psi}$ is induced from the one of the complex bounded domain,
wee see that $J_{\til{\A}_\psi}$ is integrable.
$X_n$ admits a Riemannian metric $g_{\scriptscriptstyle X_n}$ and a $2$-form $\ome_{\scriptscriptstyle X_n}$ which are given by 
$$
g_{\scriptscriptstyle X_n}(\dot\J_1, \dot\J_2)=\tr(\dot\J_1\dot\J_2)
$$
$$
\ome_{\scriptscriptstyle X_n}(\dot\J_1, \dot\J_2)=-\tr(\J \dot\J_1\dot\J_2),
$$
where $\dot\J_1, \dot\J_2\in TX_n$. The complex bounded domain $\{\, h\in \GL_n(\C)\, |\, 1_n-h^*h>0\, \}$ admits a K\"ahler structure which is given 
by 
$$
4\sqrt{-1}\pa\ol\pa \log\det (1_n-h^*h).
$$
Then under the identification $X_n\cong \{\, h\in M_n(\C)\, |\, 1_n-h^*h>0\, \}$ by using $\J_V,$ we have 
$$\ome=4\sqrt{-1}\pa\ol\pa \log\det (1_n-h^*h).$$
Then the space of maps $\wtil{\A}_\psi(M)$ inherits a Riemannian metric and a K\"ahler structure which are given by 
$$
\frac{1}{(\sqrt{-1})^{n}}\int_M \tr(\dot\J_1\dot\J_2)\lan\psi,\,\ol\psi\ran_S
$$
$$
\ome_{\wtil{\A}_\psi}=
\frac{4}{(\sqrt{-1})^{n-1}}\pa\ol\pa \int_M\log\det (1_n-h^*h) \lan\psi,\,\ol\psi\ran_S
$$
Hence $\ome_{\wtil{\A}_\psi}$ is closed.
Thus $(\wtil{\A}_\psi(M),J_{\til{\A}_\psi}, \ome_{\til{\A}_\psi}) $ is a K\"ahler manifold.
\end{proof}
Let $\wtil{\Diff}(M)$ be an extension of diffeomorphisms of $M$ by $2$-forms which is defined as 
$$
\wtil{\Diff}(M):=\{\, e^b F\,:\, F\in \Diff(M),\,\, b: 2\text{\rm -form}\, \,\}.
$$
Note that the product of $\wtil{\Diff}(M)$ is given by 
$$
(e^{b_1}F_1)( e^{b_2}F_2) :=e^{b_1+F_1^*(b_2)}F_1\circ F_2,
$$
where $F_1, F_2\in \Diff(M)$ and $b_1, b_2$ are real $2$-forms.
The action of $\wtil{\Diff}(M)$ on ${{\mathcal GC}}(M)$ by 
$$
e^{b} F_\#\circ \J\circ F_\#^{-1} e^{-b}, $$ where $F\in \Diff(M)$ acts on $\J$ by $F_\#\circ \J\circ F_\#^{-1}$ and 
 and $e^b$ is regarded as an element of SO$(\TT)$ and $F_\#$ denotes the bundle map of $\TT$ which is the lift of $F.$
We define $\wtil{\Diff(M)}_\psi$ to be a subgroup consists  of elements of $\wtil{\Diff(M)}$ which preserves $\psi$, 
$$
\wtil{\Diff}_\psi(M)=\{\, e^bF\in \wtil{\Diff}(M)\, : e^bF^*\psi=\psi\, \}.
$$
Then from (\ref{Apsi}), we have the following,
\bgn{proposition}
The symplectic structure $\ome_{\wtil{\A}_\psi}$ 
is invariant under
the action of $\psi$-preserving group $\wtil{\Diff_\psi(M)}$.\end{proposition}
We assume that type number of $\J_\psi$ is $0$, i.e., $\psi$ is given by 
$\psi=e^{b+\sqrt{-1}\ome}$, where $b$ is a real $2$-form and $\ome$ denotes a symplectic form.
We denote by Ham$_\ome(M)$ the Hamiltonian diffeomorphisms of $(M, \ome).$
\bgn{definition}
By using the $2$-form $b$, we define {\it generalized Hamiltonian diffeomorphisms} $\Ham_\ome^b(M)$ by
$$
\Ham_\ome^b(M):= \{ \, e^b F e^{-b}\, |\, F\in \Ham_\ome(M)\, \}.
$$
\end{definition}
Since $e^b F e^{-b}\psi =\psi,$ we see that $\Ham_\ome^b(M)$ is a subgroup of $\wtil{\Diff}_\psi(M).$
Thus $\Ham_\ome^b(M)$ acts on $\wtil{\A}_\psi(M)$ preserving the symplectic structure $\ome_{\wtil{A}_\psi}$
The Lie algebra of $\Ham_\ome^b(M)$ is also given by $C_0^\infty(M),$
where $C^\infty_0(M)=\{\, f\in C^\infty(M)\, |\, \int_M f\lan \psi, \,\,\ol\psi\ran_s=0\,\}.$
A Hamiltonian vector field $v$ is given by $i_v\ome=df$ for $f\in C^\infty_0(M)$. 
Then $e:=v-i_vb=\J_\psi(df)\in \TT$ is called {\it a generalized Hamiltonian element}.
Note that  we have $e\cdot\psi=\sqrt{-1}df\cdot\psi$. 

We denote by GR$_\J$ the generalized scalar curvature of $(\J, \J_\psi)$ for 
$\J \in \wtil{\A}_\psi(M),$
where GR$_\J$ is a real function on $M$.
The following is our main theorem:
\bgn{theorem}\label{th:main theorem}
There exists a moment map $\mu: {\wtil\A}_\psi(M)\to C^\infty_0(M)^*$ for the generalized Hamiltonian action which is given by the generalized scalar curvature $GR$,
$$
\lan\mu(\J), \, f\ran=(\sqrt{-1})^{-n}\int_M  (GR_\J)f \lan \psi, \,\,\ol\psi\ran_s,
$$
where $f\in C^\infty_0(M)$ and $\lan\mu(\J), \, f\ran$ denotes the coupling between $\mu(\J)$ and $f$.
 \end{theorem}
 Our proof of  Theorem \ref{th:main theorem} will be given in Section 8.
 \section{Preliminary results for proof of the main theorem}
 In order to show our main theorem, we shall rewrite the symplectic form  $\ome_{\wtil{\A}_\psi}$ by using the Clifford algebra and the pure spinors 
 $\phi_\a$ and  $\psi$. Such descriptions in terms of the Clifford algebra and pure spinors are suitable to obtain our main theorem by using Stokes' theorem.
 Let $\J$ be an almost \complex structure which is compatible with $\psi$. 
We denote by $\{\phi_\a\}$ trivializations of the canonical line bundle $K_\J$, where 
each $\phi_\a$ is a nondegenerate, pure spinor on $U_\a$ which induces the \complex structure $\J$. 
Arbitrary small deformations of almost \complex structures of $\J$ are given by the adjoint action, 
$$
e^{h(t)}\circ \J\circ e^{-h(t)},
$$
where $h(t)=h^{2,0}(t)+h^{0,2}(t)$ denotes a real section depending smoothly on a parameter $t$ which satisfies $h^{2,0}(t)\in \w^2 E_\J$ and $h^{0,2}(t)=\ol{h^{2,0}}(t)\in \w^2\ol{E_\J}$. 
Then the infinitesimal deformation $\dot{\J_h}$ is given by 
$$
\dot{\J_h}:=\frac{d}{dt}e^{h(t)}\circ \J\circ e^{-h(t)}|_{t=0}=[h, \J],
$$
where $h$ and $\J$ are regarded as elements of the Clifford algebra $\CL(\TT)$ and $[h, \J]$ denotes the commutator of $h$ and $\J$ which is identified with the bracket of Lie algebra  so$(\TT)$.
The real element $h\in \w^2E_\J\oplus \w^2\ol{E_\J}\subset\CL(\TT)$  acts on  nondegenerate pure spinors $\phi_\a$ on $U_\a$ by 
$\dot{\phi_\a}:=h\cdot\phi_\a$. 
Let $\J_1$ and $\J_2$ are two almost \complex structures which are locally induced from $\{\phi_{\a,1}\}$ and 
$\{\phi_{\a,2}\}$ respectively. 
Two real elements $h_1$ and $h_2$ give rise to  infinitesimal deformations $\dot{\J_{h_1}}$ of $\J_1$ and $\dot{\J_{h_2}}$ of $\J_2$, respectively.  We also denote by $\dot{\phi}_{\a,h_i}$ an element $h_i\cdot\phi_{\a,i}$ for $i=1,2$.


 Then the symplectic form  $\ome_{\wtil{\A}_\psi}$ as in (\ref{Apsi}) is given by
$$\ome_{\til{\A}_\psi}(\dot{\J}_{h_1}, \,\,\dot{\J}_{h_2})
=\frac{-1}{(\sqrt{-1})^n}\int_M \tr\,\J\dot{\J}_{h_1}\dot{\J}_{h_2}\lan\psi, \,\,\ol\psi\ran_s,$$
where $h_1, h_2$ are real elements of $\w^2 E_\J\oplus\w^2\ol E_\J.$
We shall begin to write the symplectic form  $\ome_{\wtil{\A}_\psi}$ in terms of 
pure spinors.
\bgn{lemma}\label{lem:saisho}
$$
\tr\,\J\dot{\J}_{h_1}\dot{\J}_{h_2}\lan\psi, \,\,\ol\psi\ran_s=-\frac{\sqrt{-1}}2\{\rho_\a^{-1}\lan \dot{\phi}_{\a, h_1},\,\,\,\ol{\dot{\phi}}_{\a, h_2}\ran_s-\rho_\a^{-1}\lan \dot{\phi}_{\a, h_2},\,\,\,\ol{\dot{\phi}}_{\a, h_1}\ran_s\}
$$
\end{lemma}
\bgn{proof}
The formula is shown by a local calculation.
Let $\{e_i\}_{i=1}^{2n}$ be a local basis of $E_\J$ such that $\lan e_i, \ol{e_j}\ran_{\scriptscriptstyle{T\oplus T^*}}=\delta_{i,j}$. 
Then the basis of $\{\ol{e_i}\}$ of $\ol{E}_\J$  is regarded as the dual basis of ${E}_\J^*$.  A real element $h\in \w^2E_\J\oplus\w^2\ol{E}_\J$ is written as 
$h=\sum_{i,j}h_{i,j}e_i\w e_j+\ol{h}_{i,j}\ol{e_i}\w\ol{e_j}$ and $\dot{\J}_{h}=[h,\J]$ is given by 
$$
\dot{\J}_{h}=h\J-\J h =\sqrt{-1}\sum_{i,j} h_{i,j} e_i\w e_j-\ol{h}_{i,j}\ol{e_i}\w\ol{e_j}
$$
Thus we have 
$$
\tr\J\dot{\J}_{h_1}\dot{\J}_{h_2}=4\sqrt{-1}\sum_{i,j} \(\ol{h}_{1,ij}h_{2,ji}-h_{1,ij}\ol{h}_{2,ji}\),
$$
where $\dot{\J}_h$ acts on $\TT$ by the adjoint $[\dot{\J}_h, \,]$.
By using the formula $\lan e\cdot\phi_\a,\,\,\ol{\phi}_\a\ran_s=-\lan \phi_\a,\,\,\,e\cdot\ol{\phi}_\a\ran_s$, we have 
$$
\lan \ol{e}_i\cdot\ol{e}_j\cdot\phi_\a,\,\,e_k\cdot e_l\cdot\ol{\phi}_\a\ran_s=\lan e_l\cdot e_k\cdot\ol{e}_i\cdot\ol{e}_j\cdot\phi_\a,\,\,\ol{\phi}_\a\ran_s.
$$
Applying $e_k\cdot\ol{e}_i+\ol{e}_i\cdot e_k=-2\lan e_k,\,\,\ol{e}_i\ran_{\scriptscriptstyle{T\oplus T^*}}$, we have 
$$
\lan \ol{e}_i\cdot\ol{e}_j\cdot\phi_\a,\,\,e_k\cdot e_l\cdot\ol{\phi}_\a\ran_s=
-4\(\delta_{kj}\delta_{li}-\delta_{ki}\delta_{lj}\)\lan \phi_\a, \,\,\,\ol\phi_\a\ran_s.
$$
Thus we obtain 
\bgn{align*}
\lan\dot{\phi}_{\a,h_1},\,\,\ol{\dot{\phi}}_{\a,h_2}\ran_s=&\sum_{i,j,k,l}\ol{h}_{1,ij}h_{2,kl}\lan \ol{e}_i\cdot\ol{e}_j\cdot\phi_\a,\,\,\,\,e_k\cdot e_l\cdot\ol{\phi}_\a\ran_s\\
=&-8\sum_{i,j}\ol{h}_{1,ij}h_{2,ji}\lan\phi_\a,\,\,\ol{\phi}_\a\ran_s
\end{align*}
We also have 
$$
\lan\dot{\phi}_{\a,h_2},\,\,\ol{\dot{\phi}}_{\a,h_1}\ran_s=-8\sum_{i,j}h_{1,ij}\ol{h}_{2,ji}\lan\phi_\a,\,\,\ol{\phi}_\a\ran_s
$$
Applying $\lan\phi_a,\,\,\,\ol{\phi}_\a\ran_s=\rho_\a\lan\psi,\,\,\,\ol{\psi}\ran_s$, we have 
$$\tr\,\J\dot{\J}_{h_1}\dot{\J}_{h_2}\lan\psi, \,\,\ol\psi\ran_s=
-\frac{\sqrt{-1}}2\{\rho_\a^{-1}\lan \dot{\phi}_{\a, h_1},\,\,\,\ol{\dot{\phi}}_{\a, h_2}\ran_s-\rho_\a^{-1}\lan \dot{\phi}_{\a, h_2},\,\,\,\ol{\dot{\phi}}_{\a, h_1}\ran_s\}
$$
\end{proof}
\bgn{proposition}\label{prop:saisho}
The symplectic form $\ome_{\wtil{\A}_\psi}$ is given by

$$
c_n^{-1}\ome_{\til{A}_\psi}(\dot{\J}_{h_1},\,\dot{\J}_{h_2})=\int_M\rho_\a^{-1}\lan \dot{\phi}_{\a, h_1},\,\,\,\ol{\dot{\phi}}_{\a, h_2}\ran_s
-\int_M\rho_\a^{-1}\lan \dot{\phi}_{\a, h_2},\,\,\,\ol{\dot{\phi}}_{\a, h_1}\ran_s
$$
where $c_n=\frac{1}{2(\sqrt{-1})^{n-1}}$ and 
$\rho_\a$ is the function as in (\ref{rho}) and $\dot{\phi}_{\a, h_i}=h\cdot\phi_{\a,i}$ for 
$i=1,2$.
\end{proposition}

\bgn{proof} {[\bf Proposition \ref{prop:saisho}]}
The result directly follows from Lemma \ref{lem:saisho} .
\end{proof}

\bgn{remark}
Since $\rho^{-1}_\a\lan\dot\phi_{\a, h_1},\ol{\dot\phi_{\a, h_2}}\ran_s=\rho^{-1}_\beta\lan\dot\phi_{\beta, h_1},\ol{\dot\phi_{\beta, h_2}}\ran_S$ for $\a,$ and $\beta$,
the $2n$-form $\rho^{-1}_\a\lan\dot\phi_{\a, h_1},\ol{\dot\phi_{\a, h_2}}\ran_s$ gives a globally defined $2n$-form on $M$. 
\end{remark}
Note that $\ome_\A$ is also written as 
\bgn{equation}c_n^{-1}\ome_\A(\dot{\J_{h_1}}, \dot{\J_{h_2}}) =
\int_M h_1\cdot\phi_\a\w\sig(\ol h_2\cdot\ol\phi_\a)\rho_\a^{-1}
-\int_Mh_2\cdot\phi_\a\w\sig(\ol h_1\cdot\ol\phi_\a)\rho_\a^{-1}
\end{equation}
\bgn{lemma}\label{lem: We have the following}
We have the following  identity with respect to $\sig$ and $d$ for a differential form $\ome$
$$d\sig \ome=
\bgn{cases}&+\sig d\a \quad (\deg\ome =\text{\rm even})\\
&-\sig d \ome\quad (\deg\ome =\text{\rm odd})
\end{cases}$$
\end{lemma}
\bgn{proof}
$\sig\ome$ is given by 
$$\sig\ome=
\bgn{cases}
+\ome,\qquad ( \deg\ome\equiv 0,1\, (\text{\rm mod} 4)\\
-\ome,\qquad ( \deg\ome\equiv 2,3\, (\text{\rm mod} 4)
\end{cases}
$$
Then the result follows.
\end{proof}

\bgn{lemma}\label{lem:e1e2}
For $e_1, e_2\in T_M\oplus T^*_M$ and differential forms $\ome_1,\ome_2\in \w^\bullet T^*_M$, we have 
$$
\lan e_1\cdot\ome_1,\,\, e_2\cdot\ome_2\ran_s+\lan e_2\cdot\ome_1,\,\,e_1\cdot\ome_2\ran_s=2\lan e_1, e_2\ran_{\scriptscriptstyle T\oplus T^*}\,\lan\ome_1,\ome_2\ran_s
$$
\end{lemma}
\bgn{proof}For $e\in \TT$ and $\ome_1,\ome_2\in\w^\bullet T^*_M$, we have 
$$
\lan e\cdot\ome_1,\,\,\ome_2\ran_S+\lan \ome_1,\,\,e\cdot\ome_2\ran_S=0.
$$
Then we have 
\bgn{align}
\lan e_1\cdot\ome_1,\,\,e_2\cdot\ome_2\ran_S+\lan e_2\cdot\ome_1,\,\,e_1\cdot\ome_2\ran_S 
=&
-\lan e_2\cdot e_1\cdot\ome_1,\,\,\ome_2\ran_S-\lan  e_1\cdot e_2\cdot\ome_1,\,\,\ome_2\ran_S\\
=&-\lan (e_2\cdot e_1+e_1\cdot e_2)\cdot\ome_1,\,\,\ome_2\ran_S\\
=&2\lan e_1, e_2\ran_{\scriptscriptstyle T\oplus T^*}\,\lan\ome_1,\ome_2\ran_S
\end{align}
\end{proof}

\bgn{lemma}\label{lem:[[h, J], t]}
Let $\theta=\theta^{1,0}+\theta^{0,1}$ be a $1$-form, where 
$\theta^{1,0}\in E_\J$ and $\theta^{0,1}\in \ol {E_\J}$. 
For $h=h^{2,0}+h^{0,2}\in \w^2E_\J\oplus \w^2\ol{E_\J},$ 
we have the following:
$$\,\,\,[[h, \J], \theta]\cdot\ol\psi =
2\sqrt{-1}\([h^{2,0}, \theta^{0,1}]-[h^{0,2}, \theta^{1,0}]\)\cdot\ol\psi
$$
\end{lemma}
\bgn{proof}
Let $\ome$ be a differential form satisfying $\J\ome=\sqrt{-1}k\ome$  for $-n\leq k\leq n$.
Then we have 
\bgn{align*}
[h, \J]\ome=&h\J\ome-\J h\ome \\
=&\sqrt{-1}k (h^{2,0}+h^{0,2})\ome-\sqrt{-1}(k-2)h^{2,0}\ome-\sqrt{-1}(k+2)h^{0,2}\ome\\
=&2\sqrt{-1}(h^{2,0}-h^{0,2})\ome
\end{align*}
Since $[h^{2,0}, \theta^{1,0}] =0$, we have 
\bgn{align*}
[[h, \J], \theta]\cdot\ol\psi=&2\sqrt{-1}[ (h^{2,0}-h^{0,2}),\,\,\theta]\cdot\ol\psi\\
=&2\sqrt{-1}\([h^{2,0}, \theta^{0,1}]-[h^{0,2}, \theta^{1,0}]\)\cdot\ol\psi
\end{align*}
\end{proof}
\bgn{remark}
If an infinitesimal deformation $[h, \J]$ preserves $\psi$, then it follows $[h, \J]\cdot\psi=0$. 
We shall consider an infinitesimal deformation $[h, \J]$ preserving $\psi$.
Since $\J\cdot\psi=0$ and $[h, \J]\cdot\psi=0$, then it follows that $[[h,\J], \theta]\cdot\psi =[h,\J]\theta\cdot\psi$. 
\end{remark}
\bgn{lemma}\label{lem:psi}
For real elements $e, \theta\in \TT$ and $h=h^{2,0}+h^{0,2}\in \w^2 E_\J\oplus \w^2\ol {E_\J}$, 
we have 
\bgn{align}\label{2lanecdotphi}
2\lan e\cdot\phi_\a, \,\,\,\theta\cdot h\cdot\ol\phi_\a\ran_s\rho_\a^{-1}
-2\lan\theta\cdot h\cdot\phi_\a,\,\,\, e\cdot\ol\phi_\a\ran_s\rho_\a^{-1}\\=
\sqrt{-1}\lan e\cdot\psi, \,\,\, [[h,\J], \theta]\cdot\ol\psi\ran_s
+\sqrt{-1}\lan \lan [[h,\J], \theta]\cdot\psi,\,\,\, e\cdot\ol\psi\ran_s
\notag
\end{align} 
\end{lemma}
\bgn{proof} 
Since $h^{2,0}\cdot\phi_\a=0$ and $[h^{2,0}, \theta^{1,0}]=0$ and $\theta^{0,1}\cdot\ol\phi_\a=0$, the left hand side of (\ref{2lanecdotphi}) is given by 
\bgn{align*}
\text{\rm (L.H.S)}=&
2\lan e\cdot\phi_\a, \,\, \theta\cdot h^{2,0}\cdot\ol\phi_\a\ran_s\rho_\a^{-1}
-2\lan\theta\cdot h^{0,2}\cdot\phi_\a,\,\, e\cdot\ol\phi_\a\ran_s\rho_\a^{-1}\\
=&2\lan e\cdot\phi_\a, \,\, \theta^{0,1}\cdot h^{2,0}\cdot\ol\phi_\a\ran_s\rho_\a^{-1}
-2\lan\theta^{1,0}\cdot h^{0,2}\cdot\phi_\a,\,\, e\cdot\ol\phi_\a\ran_s\rho_\a^{-1}\\
=&-2\lan e\cdot\phi_\a, \,\, [ h^{2,0},\theta]\cdot\ol\phi_\a\ran_s\rho_\a^{-1}
+2\lan[ h^{0,2}, \theta]\cdot\phi_\a,\,\, e\cdot\ol\phi_\a\ran_s\rho_\a^{-1}
\end{align*}
By applying Lemma \ref{lem:e1e2} and $[h^{2,0}, \theta]\cdot\phi_\a=0$, we have 
\bgn{align*}
\text{\rm (L.H.S)}=&-4\lan e, \,\, [h^{2,0}, \, \theta] \ran_{\scriptscriptstyle{T\oplus T^*}} \lan \phi_\a,\,\,\,\ol\phi_\a  \ran_s\rho_\a^{-1}
+4\lan e, \,[h^{0,2},\,\, \theta]\ran_{\scriptscriptstyle{T\oplus T^*}}\lan \phi_\a,\,\,\,\ol\phi_\a\ran_s\rho_\a^{-1}
\end{align*}
It follows from $\lan \phi_\a, \,\,\,\ol\phi_\a\ran_s=\rho_\a\lan \psi, \,\,\,\ol\psi\ran_s$ and Lemma \ref{lem:e1e2} that we have
\bgn{align*}
\text{\rm (L.H.S)}=&-2\lan e\cdot\psi, \,\,[h^{2,0}, \theta]\cdot\ol\psi\ran_s-2\lan [h^{2,0}, \theta]\cdot\psi, \,\,e\cdot\ol\psi\ran_s\\
&+2\lan [h^{0,2},\, \theta]\cdot\psi, e\cdot\ol\psi\ran_s+2\lan e\cdot\psi,\,\, [h^{0,2}, \theta]\cdot\ol\psi\ran_s
\end{align*}
From $[h^{2,0}, \J]=2\sqrt{-1}h^{2,0}$, $[h^{0,2}, \J]=-2\sqrt{-1}h^{0,2}$ , we have 
\bgn{align*}
\text{\rm (L.H.S)}=&\sqrt{-1}\lan e\cdot\psi,\,\, [[h^{2,0}, \J], \theta]\cdot\ol\psi\ran_s+\sqrt{-1}
\lan [[h^{2,0}, \J], \theta]\cdot\psi, \,\, e\cdot\ol\psi\ran_s\\
+&\sqrt{-1}\lan[[h^{0,2}, \J], \theta]\cdot\psi,\,\, e\cdot\ol\psi\ran_s+\sqrt{-1}\lan e\cdot\psi, \,\, 
[h^{0,2}, \J], \theta]\cdot\ol\psi\ran_s\\
=&\sqrt{-1}\lan e\cdot\psi,\,\,[[h,\J], \theta]\cdot\ol\psi\ran_s+
\sqrt{-1}\lan [[h, \J], \theta]\cdot\psi,\,\, e\cdot\ol\psi\ran_s
\end{align*}
Thus we have the  result.
\end{proof}
\bgn{lemma}\label{lem: N nituite}
 Let $N=N^{3,0}+N^{0,3}$ be a real section of $\w^3 E_\J\oplus \w^3\ol E_\J$, where $N^{3,0}\in \w^3 E_\J$ and $N^{0,3}=\ol{N^{3,0}}\in \w^3\ol E_J$. 
Then we have
\bgn{align*}
&\lan e\cdot\phi_\a, \,\,\,N\cdot h\cdot\ol\phi_\a\ran_s
=\lan N\cdot h\cdot\phi_\a,\,\,\,e\cdot\ol\phi_\a\ran_s=0\\
\end{align*}
\end{lemma}
\bgn{proof} Let $U_\J^k$ be the eigenspace of an eigenvalue $\sqrt{-1}k$ with respect to the action of $\J$. Then $e\cdot\phi_\a\in U^{-n+1}_\J$. Since $h^{0,2}\cdot\ol\phi_\a=0$, we have 
\bgn{align*}N\cdot h\cdot \ol\phi_\a=&(N^{3,0}\cdot h^{2,0}+N^{3,0} \cdot h^{0,2}+ N^{0,3}\cdot h^{2,0}
+N^{0,3}\cdot h^{0,2})\cdot\ol\phi_\a\\
=&N^{3,0}\cdot h^{2,0}\cdot\ol\phi_\a+ N^{0,3}\cdot h^{2,0}\cdot\ol\phi_\a
\end{align*}
Then we have $N^{3,0}\cdot h^{2,0}\cdot\ol\phi_\a \in U^{n-5}_\J$. Since $U^{n+1}_\J=\{0\}$,
$N^{0,3}\cdot h^{2,0}\cdot\ol\phi_\a=0\in U^{n+1}_\J$.
Thus we have  $\lan e\cdot\phi_\a, \,\,\,N\cdot h\cdot\ol\phi_\a\ran_s=0$.
Then it follows $\lan N\cdot h\cdot\phi_\a,\,\,\,e\cdot\ol\phi_\a\ran_s=0$. 
\end{proof}
\bgn{lemma}\label{lem: N no toriatukai}
Let $N$ be as in before. 
Then we have 
$$
N\cdot\psi=0
$$
\end{lemma}
\bgn{proof}
Let $\{e_i\}$ be a local basis of $E_\J$. 
Since $d\phi_\a=\eta_\a\cdot\phi_\a+N\cdot\phi_\a$, 
then 
we have 
\bgn{align*}
\lan N\cdot \phi_\a,\,\,\, e_i \cdot e_j\cdot e_k\cdot\ol\phi_\a\ran_s=
&\lan d\phi_\a,\,\,\, e_i \cdot e_j\cdot e_k\cdot\ol\phi_\a\ran_s\\
=&-\lan e_i\cdot e_j\cdot d\phi_\a,\,\,\, e_k\cdot\ol\phi_\a\ran_s\\
=&\lan [e_i, e_j]_{\text{cou}}\cdot\phi_\a,\,\,\, e_k\cdot\ol\phi_\a\ran_s
\end{align*}
Thus the component $N^{3,0}_{i,j,k}:=N(e_i, e_j, e_k)$ is given by 
$N^{3,0}_{i,j,k}=\lan [e_i, e_j]_{\text{cou}},\,\,\, e_k\ran_{\scriptscriptstyle{T\oplus T^*}}.$
Each $\ol e_i$ is decomposed into $\ol e_i =\ol {e}_i^+  +\ol {e}_i^-$, where $\ol e_i^\pm \in \ol E_\J^\pm$. 
From $\ol e_i^-\cdot\psi=0$ and $e_i^+\cdot\psi=0$, it suffices to show that 
$N(\ol e_i^+, \ol e_j^+, \ol e_k^+)\ol e_i^+\cdot\ol e_j^+\cdot\ol e_k^+\cdot\psi=0$ and 
$N(e_i^-,  e_j^-,  e_k^-) e_i^-\cdot e_j^-\cdot e_k^-\cdot\psi=0$.
Since $\J_\psi$ is integrable,  it follows from $[\ol e_i^+, \, \ol e_j^+]_{\text{cou}}\in \ol L_\psi$. 
From $\ol e_k^+\in \ol L_\psi$, we have
\bgn{align*}
N(\ol e_i^+, \ol e_j^+, \ol e_k^+)=\lan [\ol e_i^+, \, \ol e_j^+]_{\text{cou}},\,\ol e_k^+\ran_{\scriptscriptstyle{T\oplus T^*}}=0
\end{align*}
Thus $N(\ol e_i^+, \ol e_j^+, \ol e_k^+)\cdot\psi=0$.
We also have $N( e_i^-,  e_j^-,  e_k^-)=0$. Hence $N\cdot\psi=0$.
\end{proof}
\bgn{lemma}\label{lem: N notoriatukai2} 
If $[h, \J]\cdot\psi=0$, then we have
$$\lan e\cdot\psi,\,\,\,[[h,\J],\,\, N]\cdot\ol\psi\ran_s=
\lan [[h,\J],\,\, N]\cdot\psi,\,\,\, e\cdot\ol\psi\ran_s=0
$$
\end{lemma}
\bgn{proof}
Since $[h, \J]\cdot\psi=0$,  it follows from Lemma \ref{lem: N no toriatukai} that 
we have $[[h,\J],\,\, N]\cdot\psi=0$. Thus we have the result.
\end{proof}
\bgn{lemma}\label{lem: N notoriatukai 3}
If $[h, \J]\cdot\psi=0$, then we have
$$
\lan e\cdot N\cdot \phi_\a,\,\,\, h\cdot\ol\phi_\a\ran_s=0
$$
\end{lemma}
\bgn{proof}
Since $h=h^{2,0}+ h^{0,2}\in \w^2 E_\J\oplus \w^2\ol E_\J$ and $N=N^{3,0}+N^{0,3}\in \w^3E_\J\oplus \w^3\ol E_\J$, we have 
\bgn{align*}\lan e\cdot N\cdot\phi_\a,\,\,\, h\cdot\ol\phi_\a\ran_s
=&-\lan  N\cdot\phi_\a,\,\,\, e\cdot h\cdot\ol\phi_\a\ran_s\\
=&-\lan N^{0,3}\cdot \phi_\a,\,\,\,e^{1,0}\cdot h^{2,0}\cdot\ol\phi_\a\ran_s 
\end{align*}
Since $e^{1,0}\cdot h^{2,0}=h^{2,0}\cdot e^{1,0}$, 
we have 
\bgn{align*}
\lan e\cdot N\cdot \phi_\a,\,\,\,h\cdot\ol\phi_\a\ran_s 
=&-\lan   N^{0,3}\cdot \phi_\a,\,\,\,h^{2,0}\cdot e^{1,0}\cdot\ol\phi_\a\ran_s \\
=&\lan h^{2,0}\cdot N^{0,3}\cdot\phi_\a,\,\,\,e^{1,0}\cdot\ol\phi_\a\ran_s
\end{align*}
We denote by $[h^{2,0}, N^{0,3}]^{0,1}\in \ol E_\J$ the component  of $[h^{2,0}, N^{0,3}]$. 
Then applying Lemma \ref{lem:e1e2}, we obtain
\bgn{align*}
\lan h^{2,0}\cdot N^{0,3}\cdot\phi_\a,\,\,\,e^{1,0}\cdot\ol\phi_\a\ran_s=&
\lan [h^{2,0}\cdot N^{0,3}]^{0,1}\cdot\phi_\a,\,\,\,e^{1,0}\cdot\ol\phi_\a\ran_s\\
=&2\lan  [h^{2,0},\,\, N^{0,3}]^{0,1},\,\,\, e^{1,0}\ran_{\scriptscriptstyle{T\oplus T^*}}\lan\phi_\a,\,\,\,\ol\phi_\a\ran_s\\
=&2\lan  [h^{2,0},\,\, N^{0,3}]^{0,1},\,\,\, e^{1,0}\ran_{\scriptscriptstyle{T\oplus T^*}}\lan\psi,\,\,\,\ol\psi\ran_s\rho_\a\\
=&\lan [h^{2,0},\,\,N^{0,3}]^{0,1}\cdot\psi,\,\,\,e^{1,0}\cdot\ol\psi\ran_s\rho_\a
\\&-\lan e^{1,0}\psi,\,\,\,[h^{2,0},\,\, N^{0,3}]^{0,1}\cdot\ol\psi\ran_s\rho_\a
\end{align*}

From Lemma \ref{lem: N no toriatukai} and $h\cdot\psi=0$, we have 
$[h, N]\cdot\psi=0$. 
Thus we have $[h^{2,0},\,\, N^{0,3}]^{0,1}\cdot\ol\psi=0$ and
$[h^{2,0},\,\, N^{0,3}]^{0,1}\cdot\psi=0$. 
Hence we obtain 
$\lan e\cdot N\cdot\phi_\a,\,\,\, h\cdot\ol\phi_\a\ran_s=0$.
\end{proof}
\section{Proof of main theorem}
This section is devoted to show our main theorem: Theorem \ref{th:main theorem}.
In order to show the main theorem, it suffices to show that 
$$
\frac{d}{dt}\lan\mu(\J_t),\,\, f\ran |_{t=0}=\ome_{\wtil{\A}_\psi}(L_e\J, \dot{\J}_h),
$$
where $f$ is a generalized Hamiltonian and $e$ is a generalized Hamiltonian element satisfying
$e\cdot\psi=\sqrt{-1}df\cdot\psi$
and $\J_t$ denotes deformations of $\J$ which satisfies $\dot{\J}_h=[h, \J]$. 
A generalized Hamiltonian $f$ gives a generalized Hamiltonian element $e\in \TT$  by $e=\J_\psi df$.

Let $\{(\phi_\a, U_\a)\}$ be trivializations of 
the canonical line bundle $K_\J$, where $\{U_\a\}$ is a finite open cover of a compact manifold $M$ of dimension $2n$.
We denote by $\{\chi_\a\}$  a partition of unity such that the support of $\chi_\a$ is contained in $U_\a$.
From Proposition \ref{prop:saisho}, it is suffices to show the following: 
$$
c_n^{-1}\frac{d}{dt} \lan \mu(J_t), f\ran|_{t=0}=\int_M \lan L_e\phi_\a,\,\,\, h\cdot\ol\phi_\a\ran_s\rho_a^{-1}-\int_M \lan h\cdot\phi_\a, \,\,\, L_e\ol\phi_\a\ran_s\rho_\a^{-1}.
$$
By using the partition of unity, $f$ is given by $f=\sum_\a f_\a,$ where $f_\a=\chi_\a f$ and a generalized Hamiltonian element $e\in \TT$ is also written as $e=\sum_\a e_\a,$ where $e_\a=\J_\psi df_\a$. 
\bgn{lemma} If $U_\a\cap U_\beta\neq\emptyset$, we have 
$$\lan L_{e}\phi_\a,\,\,\, h\cdot\ol{\phi}_\a\ran_s\rho_\a^{-1}=
\lan L_{e}\phi_\beta,\,\,\, h\cdot\ol{\phi}_\beta\ran_s\rho_\beta^{-1}$$
\end{lemma}
\bgn{proof}
Since $\phi_\a=e^{\kappa_{\a,\beta}}\phi_\beta$,
the Lie derivative $L_{e}\phi_\a:= d e\cdot\phi_a+e\cdot d\phi_\a$ is given by
\bgn{align*}
L_{e}\phi_\a=&e^{\kappa_{\a,\beta}}L_e\phi_\beta+(e\cdot de^{\kappa_{\a,\beta}}+de^{\kappa_{\a,\beta}}\cdot e)\cdot\phi_\beta\\
=&e^{\kappa_{\a,\beta}}L_e\phi_\beta-2\lan  e,\,\,\, de^{\kappa_{\a,\beta}}\ran_{\scriptscriptstyle{T\oplus T^*}}\phi_\beta
\end{align*}
Since $h\in \w^2 E_\J\oplus \w^2\ol{E}_\J$, we have 
$\lan  e,\,\,\, de^{\kappa_{\a,\beta}}\ran_{\scriptscriptstyle{T\oplus T^*}}\lan \phi_\beta,\,\, h\cdot\ol{\phi}_\a\ran_s=0$. 
Since $\ol{\phi}_\a=e^{\ol{\kappa}_{\a,\beta}}\ol{\phi}_\beta$ and $\rho_\a=e^{\kappa_{\a,\beta}+\ol{\kappa}_{\a,\beta}}\rho_\beta$, 
we have the result.
\end{proof}
Thus there is a $2n$-from $F_1(e)$ such that 
$F_1(e)|_{U_\a}=\lan L_e\phi_\a,\,\,\,h\cdot\ol{\phi}_\a\ran_s\rho_\a^{-1}$. 
Since $e=\sum_\a e_\a$, it follows that 
$F_1(e)=\sum_\a F_1(e_\a)$. Since the support $e_\a$ is contained in $U_\a$, 
 we have 
$$F_1(e_\a) =\lan L_{e_\a}\phi_\a, \,\,\, h\cdot\ol{\phi}_\a\ran_s\rho_\a^{-1}.$$
Applying  Stokes' theorem and Lemma \ref{lem: We have the following}, we have 
$$
\int_M \lan d e_\a\cdot\phi_\a,\,\,\,\rho_\a^{-1}h\cdot\ol{\phi}_\a\ran_s=
\int_M  \lan e_\a\cdot\phi_\a,\,\, , d(\rho_\a^{-1}h\cdot\ol{\phi}_\a)\ran_s
$$
Thus we have 
\bgn{align*} \int_M F_1(e_\a)=&
\int_M \lan L_{e_\a}\phi_\a,\,\,\,h\cdot\ol\phi_\a\ran_s\rho_\a^{-1}\\
=&\int_M \lan d {e_\a}\cdot\phi_\a,\,\,\,h\cdot\ol\phi_\a\ran_s \rho_\a^{-1}+
\int_M \lan {e_\a}\cdot d\phi_\a,\,\,\,h\cdot\ol\phi_\a\ran_s\rho_\a^{-1}\\
=&\int_M\lan {e_\a}\cdot\phi_\a,\,\,\,dh\cdot\ol\phi_\a\ran_s\rho_\a^{-1}+
\int_M\lan{e_\a}\cdot\phi_\a,\,\,\,(d\rho_\a^{-1})\cdot h\cdot\ol\phi_\a\ran_s\\
+&\int_M\lan{e_\a}\cdot(\eta_\a+N_\a)\cdot\phi_\a,\,\,\,h\cdot\ol\phi_\a\ran_s\rho_\a^{-1}
\end{align*}
We define $F_{1-1}, F_{1-2}$ and $F_{1-3}$ by 
\bgn{align*}
F_{1-1}=&\lan {e_\a}\cdot\phi_\a,\,\,\,dh\cdot\ol\phi_\a\ran_s\rho_\a^{-1} \\
F_{1-2}=&\lan{e_\a}\cdot\phi_\a,\,\,\,(d\rho_\a^{-1})\cdot h\cdot\ol\phi_\a\ran_s \\
F_{1-3}=&\lan{e_\a}\cdot(\eta_\a+N_\a)\cdot\phi_\a,\,\,\,h\cdot\ol\phi_\a\ran_s\rho_\a^{-1}
\end{align*}

We denote by $F_2(e_\a)$ the $2n$-form $\lan h\cdot\phi_\a, \,\,\,L_{e_\a}\ol\phi_\a\ran_s\rho_\a^{-1}$. Applying Stokes' theorem again, we have 
\bgn{align*}
\int_M F_2(e_\a)=&
\int_M h\cdot\phi_\a\w\sig(L_{e_\a}\ol\phi_\a)\rho_\a^{-1}\\
=&\int_M\lan h\cdot\phi_\a,\,\, d{e_\a}\cdot\ol\phi_\a\ran_s\rho_\a^{-1}
+\int_M \lan h\cdot\phi_\a,\,\, {e_\a}\cdot d\ol\phi_\a\ran_s\rho_\a^{-1}\\
=&\int_M \lan dh\cdot\phi_\a,\,\,{e_\a}\cdot\ol\phi_\a\ran_s\rho_\a^{-1}+
\int_M \lan (d\rho_\a^{-1})\cdot h\cdot\phi_\a,\,\,{e_\a}\cdot\ol\phi_\a\ran_s\\
+&\int_M \lan h\cdot\phi_\a,\,\, {e_\a}\cdot({\eta_\a}+N)\cdot\ol\phi_\a\ran_s\rho_\a^{-1}
\end{align*}
We also define $F_{2-1}, F_{2-2}$ and $F_{2-3}$ by 
\bgn{align*}
F_{2-1}=&\lan dh\cdot\phi_\a,\,\,{e_\a}\cdot\ol\phi_\a\ran_s\rho_\a^{-1}  \\
F_{2-2}=&\lan (d\rho_\a^{-1})\cdot h\cdot\phi_\a,\,\,{e_\a}\cdot\ol\phi_\a\ran_s\\
F_{2-3}=&\lan h\cdot\phi_\a,\,\, {e_\a}\cdot({\eta_\a}+N)\cdot\ol\phi_\a\ran_s\rho_\a^{-1}
\end{align*}

 $F_1(e_\a)-F_2(e_\a)$ is divided into the following three parts 
\bgn{align*}
& F_{1-1}-F_{2-1}=\lan {e_\a}\cdot\phi_\a,\,\,\,dh\cdot\ol\phi_\a\ran_s\rho_\a^{-1}
-\lan dh\cdot\phi_\a,\,\,{e_\a}\cdot\ol\phi_\a\ran_s\rho_\a^{-1}\\
& F_{1-2}-F_{2-2}=\lan{e_\a}\cdot\phi_\a,\,\,\,(d\rho_\a^{-1}) \cdot h\cdot\ol\phi_\a\ran_s
-\lan (d\rho_\a^{-1})\cdot h\cdot\phi_\a,\,\,{e_\a}\cdot\ol\phi_\a\ran_s\\
& F_{1-3}-F_{2-3}=\lan{e_\a}\cdot(\eta_\a+N_\a)\cdot\phi_\a,\,\,\,h\cdot\ol\phi_\a\ran_s\rho_\a^{-1}
-\lan h\cdot\phi_\a,\,\, {e_\a}\cdot(\eta_\a+N_\a)\cdot\ol\phi_\a\ran_s\rho_\a^{-1}
\end{align*}
Deformations of almost \complex structures $\{\J_t\}$ are given by the action of  one parameter family $e^{ht}$ in Spin group which are induced from nondegenerate, pure spinors
$e^{ht}\cdot\phi_\a$  and we have  
\bgn{equation}\label{eht phia}
de^{ht}\cdot\phi_\a=(\eta_\a(t)+N_\a(t))\cdot e^{ht}\phi_\a,
\end{equation}
where $\eta_\a(t)\in \TT$ and $N_\a(t)\in \w^3(\TT)$ are real sections satisfying $\eta_\a(0)=\eta_\a$ and
$N_\a(0)=N_\a$. 
Taking the derivative of both sides of (\ref{eht phia}) with respect to $t$, we have 
$$
dh\cdot\phi_\a=(\dot{\eta}_\a+\dot{N}_\a)\cdot\phi_\a+(\eta_\a+N_\a)\cdot h\cdot\phi_\a,
$$
where $\dot{\eta}_\a=\frac{d}{dt}\eta_\a(t)|_{t=0}$ and $\dot{N}_\a=\frac{d}{dt}N_\a(t)|_{t=0}$.
Since the real section $\dot{{\eta_\a}}$ is decomposed into $\dot{{\eta_\a}}^{1,0}+\dot{{\eta_\a}}^{0,1}$, where 
$\dot{{\eta_\a}}^{1,0}\in E_\J$ and $\dot{{\eta_\a}}^{0,1}\in \ol E_\J$ and 
$\dot{N}$ is also decomposed in to $\sum_{p+q=3} \dot{N}^{p,q}$, where 
$\dot{N}^{p,q}\in \w^p E_\J\oplus \w^q\ol E_\J$.
Note that $\dot{N}$ is not contained in  $\w^3 E_\J\oplus \w^3\ol E_\J$ in general. 

Then we have $\dot{{\eta_\a}}\cdot\phi_\a=\dot{{\eta_\a}}^{0,1}\cdot\phi_\a$
and  $\J\phi_\a=-n\sqrt{-1}\phi_\a$. We also have $\J\dot{{\eta_\a}}\cdot\phi_\a=\J\dot{{\eta_\a}}^{0,1}\cdot\phi_\a=
(-n+1)\sqrt{-1}\dot{\eta_\a}\cdot\phi_\a.$
Then we have 
\bgn{align*}
[\J, \dot{{\eta_\a}}]\cdot\phi_\a=&\J\dot{\eta_\a}\cdot\phi_\a-\dot{\eta_\a}\J\cdot\phi_\a\\=&(-n+1)\sqrt{-1}\dot{\eta_\a}\cdot\phi_\a+n\sqrt{-1}\dot{\eta_\a}\cdot\phi_\a=\sqrt{-1}\dot{\eta_\a}\cdot\phi_\a
\end{align*}
Then we have 
\bgn{align}\label{dhcdotphi=}
dh\cdot\phi_\a=&(\dot{\eta}_\a+\dot{N}_\a)\cdot\phi_\a+(\eta_\a+N)\cdot h\cdot\phi_\a\\
=&-\sqrt{-1}[\J, \dot{\eta_\a}]\cdot\phi_\a+\dot{N}\cdot\phi_\a+({\eta_\a}+N)\cdot h\cdot\phi_\a \notag
\end{align}
We also have 
$$\dot{N}\cdot\phi_\a=-\sqrt{-1}[\J, (\dot{N}^{2,1}+\dot{N}^{1,2})]\cdot\phi_\a
-\frac13\sqrt{-1}[\J, (\dot{N}^{3,0}+\dot{N}^{0,3})]\cdot\phi_\a.$$
Since $\lan {e_\a}\cdot\phi_\a,\,\, \dot{N}\cdot\ol\phi_\a\ran_s =\lan {e_\a}\cdot\phi_\a, \,\,\, (\dot{N}^{2,1}+\dot{N}^{1,2})\cdot\ol\phi_\a\ran_s,$ we have 
\bgn{equation}\label{eq dotN}
\lan {e_\a}\cdot\phi_\a,\,\, \dot{N}\cdot\ol\phi_\a\ran_s =-\sqrt{-1}\lan {e_\a}\cdot\phi_\a,\,\, [\J, \dot{N}]\cdot\ol\phi_\a\ran_s 
\end{equation}
Substituting (\ref{dhcdotphi=}) into $F_{1-1}$ and using (\ref{eq dotN}),
we obtain
\bgn{align*}
\lan {e_\a}\cdot\phi_\a, \,dh\cdot\ol\phi_\a\ran_s\rho_\a^{-1}=&\lan {e_\a}\cdot\phi_\a,\,\,\sqrt{-1}[\J, ({\dot{{\eta_\a}}+\dot{N})]}\cdot\ol\phi_\a\ran_s\rho_\a^{-1}\\
+&\lan {e_\a}\cdot\phi_\a,\,\,{({\eta_\a}+N)\cdot h\cdot\ol\phi_\a}\ran_s\rho_\a^{-1}
\end{align*}
Thus the term $F_{1-1}$ is divided into two terms $F_{1-1-1}$ and $F_{1-1-2}$,
$$
F_{1-1}=F_{1-1-1}+F_{1-1-2}
$$
where it follows from Lemma \ref{lem: N nituite} that we have
\bgn{align*}
&F_{1-1-1}=\lan {e_\a}\cdot\phi_\a,\,\,\sqrt{-1}[\J,\,\, {\dot{{\eta_\a}}+\dot{N}]}\cdot\ol\phi_\a\ran_s\rho_\a^{-1}\\
&F_{1-1-2}=\lan {e_\a}\cdot\phi_\a,\,\,{{\eta_\a}\cdot h\cdot\ol\phi_\a}\ran_s\rho_\a^{-1}
\end{align*}
The term $F_{2-1}$ is also divided into two terms
\bgn{align*}
F_{2-1}=F_{2-1-1}+F_{2-1-2}
\end{align*}
where
\bgn{align*}
&F_{2-1-1}=\lan -\sqrt{-1}[\J, \,\,\dot{{\eta_\a}}+\dot{N}]\cdot\phi_\a,\,\,{e_\a}\cdot\ol\phi_\a\ran_s\rho_\a^{-1}\\
&F_{2-1-2}=\lan {\eta_\a}\cdot h\cdot\phi_\a,\,\,{e_\a}\cdot\ol\phi_\a\ran_s\rho_\a^{-1}
\end{align*}
By using Lemma \ref{lem:e1e2} and $\lan \phi_\a, \ol\phi_\a\ran_s=\rho_\a\lan \psi, \ol\psi\ran_s$, we obtain
\bgn{align*}
F_ {1-1-1}-F_{2-1-1}=&\sqrt{-1}\lan {e_\a}\cdot\phi_\a,\,\[\J, \,\,{\dot{{\eta_\a}}+\dot{N}]}\cdot\ol\phi_\a\ran_s\rho_\a^{-1}\\
+&\sqrt{-1}\lan [\J, \dot{{\eta_\a}}+\dot{N}]\cdot\phi_\a,\,\,{e_\a}\cdot\ol\phi_\a\ran_s\rho_\a^{-1}\\
=&2\sqrt{-1}\lan {e_\a},\,\,[\J, \dot{\eta_\a}+\dot N]\ran_{\scriptscriptstyle{T\oplus T^*}}\lan \phi_\a, \,\,\ol\phi_\a\ran_s\rho_\a^{-1}\\
=&2\sqrt{-1}\lan {e_\a},[\J, \dot{\eta_\a}+\dot N]\ran_{\scriptscriptstyle{T\oplus T^*}}\lan \psi,\,\,\ol\psi\ran_s\\
=&\lan {e_\a}\cdot\psi,\,\,\sqrt{-1}[\J, \dot{\eta_\a}+\dot N]\cdot\ol\psi\ran_s+\lan \sqrt{-1}[\J, \dot{\eta_\a}+\dot N]\cdot\psi, \,\,{e_\a}\cdot\ol\psi\ran_s
\end{align*}
It follows from Lemma \ref{lem: N no toriatukai} and $\J\cdot\psi=0$ that we have $N(t)\cdot\psi=0$. 
Thus we have $\dot{N}\cdot\psi=0$. 
It follows 
\bgn{align*}
F_ {1-1-1}-F_{2-1-1}=&\lan {e_\a}\cdot\psi,\,\,\sqrt{-1}[\J, \dot{\eta_\a}]\cdot\ol\psi\ran_s+\lan \sqrt{-1}[\J, \dot{\eta_\a}]\cdot\psi, \,\,{e_\a}\cdot\ol\psi\ran_s
\end{align*}
From Lemma \ref{lem: N notoriatukai 3}, 
the term $F_{1-3}$ is given by 
\bgn{align*}
F_{1-3}=&\lan {e_\a}\cdot({\eta_\a}+N)\cdot\phi_\a,\,\, h\cdot\ol\phi_\a\ran_s\rho_\a^{-1}=
\lan {e_\a}\cdot{\eta_\a}\cdot\phi_\a,\,\, h\cdot\ol\phi_\a\ran_s\rho_\a^{-1}\\
=&-\lan {\eta_\a}\cdot {e_\a}\cdot\phi_\a,\,\, h\cdot\ol\phi_\a\ran_s\rho_\a^{-1}
=\lan {e_\a}\cdot\phi_\a,\,\, {\eta_\a} \cdot h\cdot\ol\phi_\a\ran_s\rho_\a^{-1}=F_{1-1-2}
\end{align*}
The term $F_{2-3}$ is also given by 
\bgn{align*}
F_{2-3}=&\lan h\cdot\phi_\a, \,\, {e_\a}\cdot ({\eta_\a}+N)\cdot\ol\phi_\a\ran_s\rho_\a^{-1}=-
\lan h\cdot\phi_\a, \,\, {\eta_\a}\cdot {e_\a}\cdot\ol\phi_\a\ran_s\rho_\a^{-1}\\
=&\lan {\eta_\a}\cdot h\cdot\phi_\a, \,\, {e_\a}\cdot \ol\phi_\a\ran_s\rho_\a^{-1}=F_{2-1-2}
\end{align*}
Hence we obtain 
\bgn{align*}
F_ {\scriptscriptstyle1-1-2}+F_{\scriptscriptstyle1-3}-F_{\scriptscriptstyle2-1-2}-F_{\scriptscriptstyle2-3}=&2\lan {e_\a}\cdot\phi_\a, \,\,{\eta_\a}\cdot h\cdot\ol\phi_\a\ran_s\rho_\a^{-1}
-2\lan {\eta_\a}\cdot h\cdot\phi_\a, \,\, {e_\a}\cdot \ol\phi_\a\ran_s\rho_\a^{-1}
\end{align*}
Applying  Lemma \ref{lem:psi} and substituting $\theta={\eta_\a}$, we obtain
\bgn{align*}
F_ {\scriptscriptstyle 1-1-2}+F_{\scriptscriptstyle1-3}-F_{\scriptscriptstyle2-1-2}-F_{\scriptscriptstyle2-3}=&\sqrt{-1}\lan {e_\a}\cdot\psi, \,\, [[h, \J], {\eta_\a}]\cdot\ol\psi\ran_s
+\sqrt{-1}\lan [[h, \J], {\eta_\a}]\cdot\psi, \,\, {e_\a}\cdot\ol\psi\ran_s
\end{align*}
 We also have
 \bgn{align*}
 F_{1-2}-F_{2-2}=&\lan {e_\a}\cdot\phi_\a,\,\, d\rho_\a^{-1}\cdot h\cdot\ol\phi_\a\ran_s-\lan d\rho_\a^{-1}\cdot h\cdot\phi_\a,\,\, {e_\a}\cdot\ol\phi_\a\ran_s\\
 =&\lan {e_\a}\cdot\phi_\a,\,\, -\frac{d\rho_\a}{\rho_\a}\cdot h\cdot\ol\phi_\a\ran_s\rho_\a^{-1}-
 \lan-\frac{d\rho_\a}{\rho_\a}\cdot\ h\cdot\phi_\a,\,\, {e_\a}\cdot\ol\phi_\a\ran_s\rho_\a^{-1}
 \end{align*}
 Applying Lemma \ref{lem:psi} and substituting $\theta=\frac{d\rho_\a}{\rho_\a}$, we obtain
\bgn{align*}
F_{1-2}-F_{2-2}=&
-\frac{\sqrt{-1}}2\lan {e_\a}\cdot\psi, \,\, [[h, \J], \frac{d\rho_\a}{\rho_\a}]\cdot\ol\psi\ran_s
-\frac{\sqrt{-1}}2\lan [[h, \J], \frac{d\rho_\a}{\rho_\a}]\cdot\psi, \,\, {e_\a}\cdot\ol\psi\ran_s
\end{align*} 
Hence $F_1(e_\a)-F_2(e_\a)$ is given by the following,
\bgn{align*}F_1(e_\a)-F_2(e_\a)=
&\lan L_{e_\a}\phi_\a,\,\,\,h\cdot\ol\phi_\a \ran_s\rho_\a^{-1}-\lan h\cdot\phi_\a,\,\,\,L_{e_\a}\ol\phi_\a\ran_s\rho_\a^{-1}\\
=&\sqrt{-1}\lan {e_\a}\cdot\psi,\,\,[\J, \,\,\dot{\eta_\a}]\cdot\ol\psi\ran_s+\sqrt{-1}\lan [\J,\,\, \dot{\eta_\a}]\cdot\psi, \,\,{e_\a}\cdot\ol\psi\ran_s\\
+&\sqrt{-1}\lan {e_\a}\cdot\psi, \,\, [[h, \J], {\eta_\a}]\cdot\ol\psi\ran_s
+\sqrt{-1}\lan [[h, \J], {\eta_\a}]\cdot\psi, \,\, {e_\a}\cdot\ol\psi\ran_s\\
-&\frac{\sqrt{-1}}2\lan {e_\a}\cdot\psi, \,\, [[h, \J], \frac{d\rho_\a}{\rho_\a}]\cdot\ol\psi\ran_s
-\frac{\sqrt{-1}}2\lan [[h, \J], \frac{d\rho_\a}{\rho_\a}]\cdot\psi, \,\, {e_\a}\cdot\ol\psi\ran_s
\end{align*}
Since ${e_\a}$ is a generalized Hamiltonian element satisfying
${e_\a}\cdot\psi=\sqrt{-1}d{f_\a}\cdot\psi$, we have 
\bgn{align*}
F_1(e_\a)-F_2(e_\a)=
&-\lan d{f_\a}\cdot\psi,\,\,[\J, \dot{\eta_\a}]\cdot\ol\psi\ran_s+\lan [\J, \dot{\eta_\a}]\cdot\psi, \,\,d{f_\a}\cdot\ol\psi\ran_s\\
-&\lan d{f_\a}\cdot\psi, \,\, [[h, \J], {\eta_\a}]\cdot\ol\psi\ran_s
+\lan [[h, \J], {\eta_\a}\cdot\psi, \,\, d{f_\a}\cdot\ol\psi\ran_s\\
+&\frac12\lan d{f_\a}\cdot\psi, \,\, [[h, \J], \frac{d\rho_\a}{\rho_\a}]\cdot\ol\psi\ran_s
-\frac12\lan [[h, \J], \frac{d\rho_\a}{\rho_\a}]\cdot\psi, \,\, d{f_\a}\cdot\ol\psi\ran_s
\end{align*}

The action of Spin group preserves the form $\lan\,\,,\,\,\ran_s$.
Since deformations  $J_t:=e^{h(t)}\circ \J\circ e^{-h(t)}$ is given by the action of Spin group $e^{h(t)}$,
thus $\rho_\a$  does not depend on $t$. Recall $\dot{\J}=[h,\J]$. 
Then $F_1(e_\a)-F_2(e_\a)$ is given by the following derivative at $t=0$,
\bgn{align*}
F_1(e_\a)-F_2(e_\a)=&-\frac{d}{dt}\lan d{f_\a}\cdot\psi,\,\, [\J_t,\,\,{\eta_\a}(t)]\cdot\ol\psi\ran_s
+\frac{d}{dt}\lan [\J_t, \,\,{\eta_\a}(t)]\cdot\psi,\,\,d{f_\a}\cdot\ol\psi\ran_s\\
+&\frac12\frac{d}{dt}\lan d{f_\a}\cdot\psi,\,\, [\J_t,\frac{d\rho_\a}{\rho_\a}]\cdot\ol\psi\ran_s-\frac12\frac{d}{dt}
\lan [\J_t, \frac{d\rho_\a}{\rho_\a}]\cdot\psi,\,\,d{f_\a}\cdot\ol\psi\ran_s
\end{align*}
Since we consider deformations preserving $\psi$,
we have $\J_t\cdot\psi=0$. Thus we have $ [\J_t,\,\,{\eta_\a}(t)]\cdot\ol\psi=\J_t{\eta_\a}(t)\cdot\ol\psi$.
The support of $f_\a$ is contained in $U_\a$.
Applying Stokes' theorem, we obtain

\bgn{align*}
\int_M F_1(e_\a)-\int_MF_2(e_\a)
=&\frac{d}{dt}\int_M \lan f_\a\psi, \,\,\, d(-\J_t {\eta_\a}(t)+\frac12\J_td\log\rho_\a)\cdot\ol\psi\ran\\
+&\frac{d}{dt}\int_M \lan d(\J_t{\eta_\a}(t)-\frac12 \J_t\log\rho_\a)\cdot\psi, \,\,\,f_\a\ol\psi\ran_s
\end{align*}
Since 
$d(\J\eta_\a+\frac12\J d\log\rho_\a)\cdot\ol{\psi}$ is a globally defined $d$-closed $2n$-form,
we have 
\bgn{align*}
\int_M F_1(e)-F_2(e)=&
\sum_\a\int_M F_1(e_\a)-\sum_\a\int_MF_2(e_\a)\\
=&\sum_\a\frac{d}{dt}\int_M \lan f_\a\psi, \,\,\, (-d\J_t {\eta_\a}(t)+\frac12d\J_td\log\rho_\a)\cdot\ol\psi\ran\\
+&\sum_\a\frac{d}{dt}\int_M \lan( d\J_t{\eta_\a}(t)-\frac12 d\J_t\log\rho_\a)\cdot\psi, \,\,\,f_\a\ol\psi\ran_s\\
=&\frac{d}{dt}\int_M \lan f\psi, \,\,\, (-d\J_t {\eta_\a}(t)+\frac12d\J_td\log\rho_\a)\cdot\ol\psi\ran\\
+&\frac{d}{dt}\int_M \lan( d\J_t{\eta_\a}(t)-\frac12 d\J_t\log\rho_\a)\cdot\psi, \,\,\,f\ol\psi\ran_s
\end{align*}
Hence we obtain 
\bgn{align*}
c_n^{-1}\ome_{\wtil{\A}_\psi}(L_e\J, \,\dot{\J}_h) =&\frac{d}{dt}\int_M \lan f\psi,\,\,\,  d(-\J_t{\eta_\a}(t)+\frac12\J _td\log\rho_\a)\cdot\ol{\psi}\ran_s\\
&+\frac{d}{dt}\int_M\lan d(\J_t{\eta_\a}(t)-\frac12\J_t\log\rho_\a)\cdot\psi,\,\,\, f\ol\psi\ran_s
\end{align*}
Thus it follows from (\ref {GR in general})  that the moment map $\mu$ is given by
\bgn{align*}
\lan \mu(\J), f\ran =&c_n\int_M \lan f\psi,\,\,\, d(-\J{\eta_\a}+\frac12\J d\log\rho_\a)\cdot\ol\psi\ran_s\\
+&c_n\int_M\lan d(\J{\eta_\a}-\frac12\J\log\rho_\a)\cdot\psi, \,\,\, f\ol\psi\ran_s\\
=&(\sqrt{-1})^{-n}\int_M f (GR_\J) \lan \psi, \,\,\ol\psi\ran_s
\end{align*}
Hence we obtain the result.

\section{Deformations of generalized K\"ahler structures with constant generalized scalar curvature }
\bgn{definition}
If the generalized scalar curvature GR of a generalized K\"ahler structure
$(\J, \J_\psi)$ is constant,  then $(\J, \J_\psi)$ is called a generalized K\"ahler structure with constant generalized scalar curvature, that is, 
$$\text{\rm GR}=\lam\,\, ( \text{\rm constant}),$$ where 
$\dstyle{\lam =n\frac{ c_1(K_\J)\cup [\ome]^{n-1}}{[\ome^n]}}$
\end{definition}
\bgn{theorem}
Infinitesimal deformations of generalized K\"ahler structures with constant generalized scalar curvature are given by 
the cohomology group 
$$\ker\ol\pa_\J\cap (\ol{E}_\J^+{\scriptstyle\w}\ol{E}_\J^-)
 /\text{\rm Im} \ol\pa_\J^+\ol\pa_\J^-\cap (\ol{E}_\J^+{\scriptstyle\w}\ol{E}_\J^-)
$$ of the following elliptic complex :
$$
0\to C^\infty_\C(M)\overset{\ol\pa_\J^+\ol\pa_\J^-}{\longrightarrow}\ol{E}_\J^+{\scriptstyle\w}\ol{E}_\J^-
\overset{\ol\pa_{\scriptstyle \J}}{\longrightarrow}(\w^2\ol{E}_\J^+{\scriptstyle\w}\ol{E}_\J^-)\oplus(\ol{E}_\J^+{\scriptstyle\w}\w^2\ol{E}_\J^-)
\overset{\ol\pa_{\scriptstyle \J}}{\longrightarrow}\cdots
$$
Since the cohomology group is finite dimensional for a compact manifold $M$, Infinitesimal deformations 
are also finite dimensional.
\end{theorem}
\bgn{proof}
Let ${\A}_\psi(M)$ be the set of \complex structures which are compatible with $\psi$ and 
Ham$_\psi(M)$ the generalized Hamiltonian group which acts on $\A_\psi(M)$. 
The orbit of Ham$_\psi(M)$ on $\A_\psi(M)$ is denoted by ${\mathcal O}_{Ham}(M)$.
Let $\J\in \A_\psi(M)$ be a \complex structure such that $(\J, \J_\psi)$ admits constant generalized scalar curvature. 
The formal tangent space of $\A_\psi(M)$ at $\J$ is given by $\e\in \ol E_\J^+\w \ol E_\J^-$ satisfying $\ol\pa_\J \e=0$,
since deformations of $\J$ preserves $\psi$. 
Since generalized K\"ahler structures with constant generalized scalar curvature are given by 
the inverse image $\mu^{-1}(0)$ of the moment map $\mu$ for the action of  Ham$_\psi(M)$, 
infinitesimal deformations are the orthogonal complement of the direct sum of $(T_\J{\mathcal O}_{Ham}\oplus 
\J T_\J{\mathcal O}_{Ham})$, where $T_\J{\mathcal O}_{Ham}$ denotes the tangent space of the orbit ${\mathcal O}_{Ham}$ at $\J$.
The $T_\J{\mathcal O}_{Ham}$ consists of $L_e\J$ for Hamiltonian element $e=e^{1,0}+e^{0,1}\in \TT$, 
where $e^{1,0}\in E_\J$ and $e^{0,1}\in \ol E_\J$. 
Thus $T_\J{\mathcal O}_{Ham}$ is given by $\{\ol\pa_\J e^{0,1}\, |\, e: \text{\rm Hamiltonian element}\}$.
Since a Hamiltonian element $e$ is given by 
$e=\J_\psi df$ for a hamiltonian $f$, we have 
\bgn{align*}
 \ol\pa_\J e^{0,1}=&\ol\pa_\J (\J_\psi df)^{0,1} =
 \sqrt{-1}(\ol\pa_\J^++ \ol\pa_\J^-)(\ol\pa_\J^+ - \ol\pa_\J^-)f \\
 =&-2\sqrt{-1}\,\,\ol\pa_\J^+\ol\pa_\J^- f
\end{align*}
Since $\J$ acts on $2\sqrt{-1}\,\ol\pa_\J^+\ol\pa_\J^- f\in \w^2\ol E_\J$ by the multiplication of 
$2\sqrt{-1}.$
Thus we have the complexification,
$$
T_\J{\mathcal O}_{Ham}\oplus \J T_\J{\mathcal O}_{Han}=\{\, -2\sqrt{-1}\,\,\ol\pa_\J^+\ol\pa_\J^- F\, |\, F \in C^\infty_\C(M)\, \}
$$
Hence infinitesimal deformations of generalized K\"ahler structures with constant generalized scalar curvature are given by 
the cohomology group 
$$\ker\ol\pa_\J\cap (\ol{E}_\J^+{\scriptstyle\w}\ol{E}_\J^-)
 /\text{\rm Im} \ol\pa_\J^+\ol\pa_\J^-\cap (\ol{E}_\J^+{\scriptstyle\w}\ol{E}_\J^-)
$$
The ellipticity of the complex follows from checking its symbol complex.
Hence we obtain the result.
\end{proof}
\bgn{example}
Let $S$ be a K3 surface and $(\J_J, \J_\psi)$ a generalized K\"ahler structure induced from 
a Ricci flat K\"ahler structure. 
We have the generalized Hodge decomposition $H^\bullet (S) =\oplus H^{p,q}$,  
\bgn{center}
\xymatrix{&&H^{0,2}&&\\
&H^{-1,1}&&H^{1,1}&\\
H^{-2,0}& & H^{0,0} && H^{2,0}\\
&H^{-1,-1}&&H^{1,-1}&\\
&&H^{0,-2}&&}
\end{center}
Then infinitesimal deformations of generalized K\"ahler structures with vanishing generalized scalar curvature are given by 
$H^{0,0}(S),$ where $\dim H^{0,0}=20.$ 
In the cases of ordinary K3 surfaces, deformations of complex structures with vanishing Ricci tensor preserving 
a symplectic structure is 19 dimensional. Hence there is one more dimensional deformations which deform to 
generalized K\"ahler structures of type $(0,0)$ which is discussed next section.
\end{example}

\section{Generalized K\"ahler structures of type $(0,0)$}
\bgn{definition}
A generalized K\"ahler structure of type $(0,0)$ is a generalized K\"ahler structure $(\J_\phi, \J_\psi)$
which is induced from a pair 
$$(\phi=e^{B+\sqrt{-1}\ome_1}, \psi=e^{\sqrt{-1}\ome_2})$$ which consists of $d$-closed, nondegenerate, pure spinors of symplectic types, where $B$ is a real $d$-closed $2$-form and both $\ome_1$ and $\ome_2$ are real symplectic forms, respectively. 
\end{definition}
\bgn{proposition}
A pair $(\phi=e^{B+\sqrt{-1}\ome_1}, \psi=e^{\sqrt{-1}\ome_2})$ gives a generalized K\"ahler structure if and only if 
$(\phi, \psi)$ satisfies the followings :\\
(1) $\ome_\C^{\pm} :=B+\sqrt{-1}(\ome_1\mp\ome_2)$ defines complex structures $I_\pm$ such that $\ome_\C^\pm$ are 
$d$-closed holomorphic symplectic forms with respect to $I_\pm$ respectively.\\
(2) $\ome_2$ is tame w.r.t both $I_\pm$.
\end{proposition}
\bgn{proof}Let $E_\phi$ be the eigenspace with eigenvalue $-\sqrt{-1}$ with respect to $\J_\phi$ and 
$\ol E_\phi$ the complex conjugate of $E_\phi$.
We denote by $E_\psi$ the eigenspace with eigenvalue $-\sqrt{-1}$ with respect to $\J_\psi$ and $
\ol E_\psi$ is the complex conjugate of $E_\psi$.
Then we have 
\bgn{align*}
& E_\phi=\{\ v-i_v(B+\sqrt{-1}\ome_1)\, |\, v\in T_M^\C\, \}, \qquad E_\psi =\{\, u-\sqrt{-1}i_u\ome_2\,|\, u\in T_M^\C\,\}\\
&\ol E_\phi=\{\ v-i_v(B-\sqrt{-1}\ome_1)\, |\, v\in T_M^\C\, \}, \qquad \ol E_\psi =\{\, u+\sqrt{-1}i_u\ome_2\,|\, u\in T_M^\C\}
\end{align*}
The condition $\J_\phi \J_\psi=\J_\psi\J_\phi$ is equivalent to the followings : 
$$\dim_\C  E_\phi\cap E_\psi =\dim_\C E_\phi\cap \ol E_\psi=n.$$ 
Thus $u-\sqrt{-1}i_u\ome_2\in E_\phi\cap E_\psi$ if and only if 
$u-\sqrt{-1}i_u\ome_2=u-i_u(B+\sqrt{-1}\ome_1)$. 
Hence $u-\sqrt{-1}i_u\ome_2\in E_\phi\cap E_\psi$ if and only if $i_u(B+\sqrt{-1}(\ome_1-\ome_2))=0$. 
Thus $\ker (B+\sqrt{-1}(\ome_1-\ome_2)):=\{\, u\in T_M^\C\, |\, i_u(B+\sqrt{-1}(\ome_1-\ome_2)=0\, \}$ is 
$n$ dimensional if and only if $\dim E_\phi\cap E_\psi=n$.
If $u\in E_\phi$, then it follows from $E_\phi\cap \ol E_\phi=\{0\}$ that we have $u\neq \ol u$. 
Thus we see that 
$$\ker(B+\sqrt{-1}(\ome_1-\ome_2))\cap \ol {\ker (B+\sqrt{-1}(\ome_1-\ome_2))}=\{0\}.$$
Hence $\ome_\C^+:=B+\sqrt{-1}(\ome_1-\ome_2)$ defines a complex structure $I_+$ such that $\ome_\C^+$ is a 
holomorphic symplectic form with respect to $I_+$.
We also see that $\ker (B+\sqrt{-1}(\ome_1+\ome_2)):=\{\, u\in T_M^\C\, |\, i_u(B+\sqrt{-1}(\ome_1+\ome_2))=0\, \}$ is 
$2n$ dimensional if and only if $\dim E_\phi\cap \ol E_\psi=n$.
Thus $\ome_\C^-:=B+\sqrt{-1}(\ome_1+\ome_2)$ defines a complex structure $I_-$ such that 
$\ome_\C^-$ is a holomorphic symplectic form with respect to $I_-$.
Hence the condition $[\J_\phi, \J_\psi]=0$ is equivalent to the condition (1). 
The eigenspace with eigenvalue $\pm 1$ with respect to $G:=\J_\phi\J_\psi$ are denoted by $C_\pm$, respectively. 
Then we have $C_+^\C=(E_\phi\cap E_\psi)\oplus (\ol E_\phi\cap \ol E_\psi)$ and 
$C_-^\C=(E_\phi\cap \ol E_\psi )\oplus (\ol E_\phi\cap E_\psi)$. 
For $ u\in \ker\ome_\C^+=T^{0,1}_{I_+}$,we have 
\bgn{align*}
G( u-\sqrt{-1}i_u\ome_2,\, \ol u+\sqrt{-1}i_{\ol u}\ome_2,\,)=&\lan  u-\sqrt{-1}i_u\ome_2, \, \ol u+\sqrt{-1}i_{\ol u}\ome_2,\ran\\
=&-2\sqrt{-1}\ome_2(u, \ol u)
\end{align*}
For $u\in \ker\ome_\C^-=T^{0,1}_{I_-}$, we also have 
\bgn{align*}
G(u+\sqrt{-1}i_u\ome_2, \, \ol u-\sqrt{-1}i_{\ol u}\ome_2)=&-\lan u+\sqrt{-1}i_u\ome_2, \, \ol u-\sqrt{-1}i_{\ol u}\ome_2\ran\\ =&-2\sqrt{-1}\ome_2(u, \ol u)
\end{align*}
Thus $G=\J_\phi\J_\psi$ gives a generalized metric if and only if $-\sqrt{-1}\ome_2(u, \ol u)> 0$ for all $u\neq 0\in T^{0,1}_{I_\pm}$.
A symplectic structure is tame with respect to $I_\pm$ if and only if 
$\ome_2(x, I_\pm x)> 0$ for  every real tangent $x\neq 0\in T_M$. 
Since $-\sqrt{-1}\ome_2( x-\sqrt{-1}I_\pm x, \, x+\sqrt{-1}I_\pm x)=2\ome_2(x, I_\pm x)$, 
Hence $G:=\J_\phi\J_\psi$ gives a generalized metric if and only if $\ome_2$ is tame with respect to $I_\pm$.
Hence we obtain the result.
\end{proof}
\bgn{remark}
On a $4$ dimensional manifold, the condition (1) is equivalent to the followings 
$$
B\w \ome_1=B\w\ome_2=\ome_1\w\ome_2=0, \quad B\w B=\ome_1\w\ome_1+\ome_2\w\ome_2\neq0.
$$ 
\end{remark}
In the case of a generalized K\"ahler structure of type $(0,0),$
the GRic and GR are explicitly written. 
\bgn{proposition}
For a generalized K\"ahler structure of type $(0,0)$, GRic and GR are given by 
\bgn{align*}
\text{GRic}=&-d B\ome_1^{-1}(d\log \frac{\ome_1^n}{\ome_2^n})\\
(\text{GR})\ome_2^n=&n\ome_2^{n-1}\w d B\ome_1^{-1}(d\log \frac{\ome_1^n}{\ome_2^n}),
\end{align*}
where $B: T_M\to T^*_M$ and $\ome_i^{-1}: T^*_M\to T_M$ ($i=1,2$) and the composition $B\ome_1^{-1}$
is an endomorphism of $T^*_M$ and then $B\ome_1^{-1}(d\log \frac{\ome_1^n}{\ome_2^n})$ is a 
$1$-form and then the exterior derivative of $B\ome_1^{-1}(d\log \frac{\ome_1^n}{\ome_2^n})$  is a $2$-form which is
the GRic form.
\end{proposition}
\bgn{proof}
Substituting $\rho_\a=\frac{\ome_1^n}{\ome_2^n}$ and $\eta_\a=0$ into (\ref{GRic}), we have the result.
\end{proof}

\bgn{example}[HyperK\"ahler str.]
Let $(g, I,J,K)$ be a hyperK\"ahler structure with three K\"ahler forms $(\ome_I, \ome_J, \ome_K)$.
We define $B$ and two symplectic forms $\ome_1, \ome_2$ by 
$$B=\ome_J,\,\,\,\, \ome_1=\frac12(\ome_I+\ome_K),\,\,\,\, \ome_2=\frac12(\ome_I-\ome_K).$$
Then 
$(\phi=e^{B+\sqrt{-1}\ome_1},\,\,\psi=e^{\sqrt{-1}\ome_2})$ is a generalized K\"ahler structure which satisfies  GRic$=0$.
\end{example}

\section{Generalized K\"ahler-Einstein structures}
\bgn{definition}A generalized K\"ahler structure $(\J_\beta, \psi=e^{b+\sqrt{-1}\ome})$ is a generalized K\"ahler-Einstein if 
we have the following:
$$
\text{\rm GRic} =\lam\ome $$
for constants $\lam$.
\end{definition}
In the case of generalized K\"ahler structure of type $(0,0)$, the generalized K\"ahler-Einstein condition implies that 
$\ome_1^n=\ome_2^n,$ where the Einstein constant is zero.
\section{Generalized K\"ahler-Einstein structures constructed from holomorphic Poisson deformations}
Let $(M, J, \ome)$ be a K\"ahler manifold with an ordinary complex structure $J$ and a K\"ahler structure $\ome$. 
We assume that the $m$-dimensional torus $T$ acts on $M$ preserving the K\"ahler structure $(J, \ome)$ on $M$ and 
there exists a moment map $\mu_T: M \to (t^m)^*$ for the action of $T$, where we assume $m\geq 2$. 
Let $\{\xi_i\}_{i=1}^m$ be a basis of the Lie algebra $t^m$ of the Torus $T$ and 
$\{V_i\}_{i=1}^m$ the corresponding real vector fields which are generated by $\{\xi_i\}_{i=1}^m$. 
Each $V_i$ is decomposed into $V_i^{1,0}+V^{0,1}_i$, where $V_i^{1,0}\in T^{1,0}_J$ and 
$V^{0,1}_J\in T^{0,1}_J$. Since $\{V_i\}_{i=1}^m$ are commuting vector fields, we have a real Poisson structure $\beta_\R$ by 
$$
\beta_\R=\sum_{i,j} \lam_{i,j}V_i\w V_j,
$$  
where $\lam_{i,j}$ are constants.
Holomorphic vector fields $\{V_i^{1,0}\}_{i=1}^m$ also gives a holomorphic Poisson structure 
$\beta=\sum\lam_{i,j}V_i^{1,0}\w V_j^{1,0}$.
Let $(\J_J, J_\psi)$ be the generalized K\"ahler structure coming from the ordinary K\"ahler structure 
$(J, \ome)$, where $\psi =e^{\sqrt{-1}\ome}$. 
Let $\{\phi_\a\}$ be trivializations of $K_{\J_J}$, that is, each $\phi_\a$ is a holomorphic $n$-form with respect to $J$. 
Then the action of $e^{\beta_\R}$ on each $\phi_\a$ coincides with the action of $e^{\beta}$ on $\phi_\a$, that is, 
$$ e^{\beta_\R}\cdot \phi_a=e^{\beta}\cdot\phi_\a.$$
Thus the action of $e^{\beta_\R}$ on $\J_J$ gives Poisson deformations of $\J_{\beta t}$. 
Then the action of $e^{\beta_\R}$ gives deformations of almost generalized K\"ahler structures 
$$(\J_{\beta t}, \J_{\psi_t}):=(e^{\beta_{\R} t}\J_J e^{-\beta_{\R} t}, \,\, e^{\beta_{\R }t} \J_\psi e^{-\beta_{\R} t}),$$ where $\J_{\psi_t}$ are almost \complex structures induced from $\psi_t=e^{\beta_\R t}\cdot\psi$.
\bgn{theorem}\label{torus GK} 
Let $\mu_{T, i}$ be the function which is the coupling $\lan \mu_T, \xi_i\ran$ of  the moment map $\mu_T$ and $\xi_i \in t^m$. Then 
$\psi_t$ is given by 
$$
\psi_t =exp{(-\sum_{i,j} \lam _{i,j} \,d\mu_{T, i}\w d\mu_{T, j}+ \sqrt{-1}\,\ome)}
$$
Thus
$d\psi_t=0$ and $(\J_{\beta t}, \J_{\psi_t})$ are deformations of generalized K\"ahler structures.
\end{theorem}
\bgn{proof}
The exponential $e^{\beta_\R}$ is given by 
$e^{\beta_\R}=\prod_{i,j} e^{\lam_{i,j}V_i\w V_j}=\prod_{i,j} (1+\lam_{i,j}V_i\w V_j).$

Since $\ome (V_i, V_j)=0$ and $i_{V_i}\ome=d\mu_{T, i}$, we have 
$$
V_i\w V_j \cdot\psi =-d\mu_{T,i}\w d\mu_{T, j}\w e^{\sqrt{-1}\ome}
$$
Since $i_{V_i}d\mu_{T, j}=0$, we have 
\bgn{align*}
e^{\beta_\R}\cdot\psi=&\prod_{i,j}(1+\lam_{i,j}V_i\w V_j)\cdot\psi=
\prod_{i,j}(1-\lam_{i,j}d\mu_{T,i}\w d\mu_{T, j})\cdot\psi\\
=&\prod_{i,j}e^{-\lam_{i,j}d\mu_{T,i}\w d\mu_{T, j}}\cdot\psi\\
=&exp{(-\sum_{i,j}\, \lam _{i,j} d\mu_{T, i}\w d\mu_{T, j}+ \sqrt{-1}\,\ome)}.
\end{align*}
Thus $\psi_t$ is $d$-closed and $\J_{\psi_t}$ are \complex structures. 
Thus we have the result.
\end{proof}

\bgn{proposition}\label{prop: GKE}
Let $(X, J,\ome)$ be a K\"ahler-Einstein manifold which admits an action of real torus 
$T^m$ $(m\geq 2)$ preserving the K\"ahler structure $(J, \ome)$. We assume that there exists a moment map for the action of $T^m$.
We denote by $\{\xi_i\}_{i=1}^m$ a basis of the Lie algebra $t^m$ which yields vector fields $\{V_i\}_{i=1}^m$.
We assume that  $\beta_\R:=\sum_{i,j}\lam_{i,j}V_i\w V_j$ is a nontrivial real Poisson structure for some constants $\lam_{i,j}$. 
Then there exist nontrivial deformations of generalized K\"ahler-Einstein manifolds
$(\J_{\beta t}, \psi_t),$ where $\{\J_{\beta t}\}$ are Poisson deformations of $\J_J,$ where 
$\beta$ is the holomorphic Poisson structure given by 
$$
\beta=\sum_{i,j}\lam_{i,j} V^{1,0}_i\w V^{1,0}_j
$$ 
and $V_i=V^{1,0}+ V^{0,1}$ and  $V^{1,0}_i\in T^{1,0}_J,$ and $ V^{0,1}\in T^{0,1}_J.$
\end{proposition}
\bgn{proof}

It suffices to show Proposition \ref{prop: GKE}  in the case of $\beta_\R=V_1\w V_2$ which is a real Poisson structure given by the wedge of $V_1$ and $V_2$. 
Let $\{\phi_a\}$ be  trivializations of the canonical line bundle $K_J$ which are given by the ordinary holomorphic $n$-forms. The action of $T^m$ preserves the complex structure $J$ and the canonical line bundle.
Thus the action of $V_1$ and $V_2$ are the representations of weights $n_1$ and $n_2$, respectively, that is,
$$
L_{V_1}\phi_\a=\sqrt{-1}n_1\phi_a\, \quad L_{V_2}\phi_\a=\sqrt{-1}n_2\phi_\a
$$
From $[V_1, V_2]=0$ and $d\phi_\a=0$, it follows that we have
\bgn{align*}
d(\beta_\R\cdot\phi_\a)=& d(V_1\w V_2)\cdot\phi_\a=L_{V_1} V_2\cdot\phi_\a-V_1 d V_2\cdot\phi_\a\\
=&V_2\cdot L_{V_1}\phi_\a-V_1\cdot L_{V_2}\phi_\a \qquad\qquad \\
=&\sqrt{-1}n_1 V_2\cdot\phi_\a-\sqrt{-1}n_2 V_1\cdot\phi_\a
\end{align*}
Since $V_i\cdot \beta_\R=V_i\cdot V_1\cdot V_2=0$, we have
\bgn{align*}
de^{\beta_\R}\phi_\a=&(\sqrt{-1}n_1V_2-\sqrt{-1}n_2V_1)\cdot e^{\beta_\R}\phi_\a
\end{align*}
Since $(V_1-\sqrt{-1}\J_\beta V_1)\cdot e^{\beta_\R}\phi_\a=(V_2-\sqrt{-1}\J_\beta V_2)\cdot e^{\beta_\R}\phi_\a=0$, we have 
$$
de^{\beta_\R}\phi_\a=(-n_1\J_\beta V_2+n_2\J_\beta V_1)\cdot e^{\beta_\R}\phi_\a.
$$
Since $\J_\beta V_i =e^{\beta} \J_J e^{-\beta} V_i =\J_J V_i = J V_i$, we also have 
$$
de^{\beta_\R}\phi_\a=(-n_1J V_2+n_2J V_1)\cdot e^{\beta_\R}\phi_\a.
$$
Since $(-n_1JV_2+n_2JV_1)$ is a real section, it follows that $\eta_\a=(-n_1JV_2+n_2JV_1) $
$\,{}^{*3}$
\footnote{${}^{*3}$  
If $\beta_\R=\sum_{i,j}\lam_{i,j}V_i\w V_j$, then $\eta_\a =\sum_{i,j} \lam_{i,j}(-n_jJV_i+n_iJ V_j)$ }
Let $\mu_T$ be the moment map for the action of $T^m$. Then $\mu_{T,i}$ is denoted by $\lan\mu_T, \xi_i\ran$.
Since $\ome(V_1, V_2)=0$, $i_{V_i}\ome=d\mu_i$, we have 
\bgn{align}
e^{\beta_\R}\cdot\psi=&e^{\beta_\R}\cdot e^{\sqrt{-1}\ome}=e^{\sqrt{-1}\ome}+V_1\w V_2\cdot e^{\sqrt{-1}\ome}\\
=&e^{\sqrt{-1}\ome}-d\mu_{T,1}\w d\mu_{T,2}w e^{\sqrt{-1}\ome}=\exp{(-d\mu_{T,1}\w d\mu_{T,2}+\sqrt{-1}\ome)}\notag
\end{align}
Since $\psi_{\beta_\R}:=e^{\beta_\R}\cdot\psi$ is $d$-closed, then $\J_\psi$ is integrable. 
Hence $(\J_{\beta_\R},\, \J_\psi)$ is a generalized K\"ahler structure.
Since $\eta_\a\in T_M$, it follows that $e^{\beta_\R}\eta_\a e^{-\beta_\R}=\eta_a$. Thus we have 
\bgn{align*}
d\J_{\beta_\R} \eta_\a\cdot\psi_{\beta_\R}=&de^{\beta_\R}\J_Je^{-\beta_\R}\eta_\a e^{\beta_\R}\cdot\psi\\
=&de^{\beta_\R}J \eta_\a \psi=de^{\beta_\R}(n_1V_2- n_2V_1)e^{-\beta_\R}e^{\beta_\R}\psi\\
=&d(n_1V_2-n_2V_1)e^{\beta_\R}\cdot\psi
\end{align*}
Since $V_i\cdot d\mu_j =0$ for $i, j=1,2$, we have
$dV_ie^{\beta_\R}\cdot\psi =dV_i\psi$. 
Thus we have 
\bgn{equation}
d\J_{\beta_\R} \eta_\a\cdot\psi_{\beta_\R}=d(n_1V_2-n_2V_1)\cdot\psi=\sqrt{-1}d (n_1 i_{V_2}\ome - n_2i_{V_1}\ome)\cdot\psi
=0.
\end{equation}

We calculate the term $d\J_{\beta_\R}d\log\rho_\a\cdot\ol\psi_{\beta_\R}.$
Since $V_i$ preserves the function $\rho_\a$, we have $L_{V_i}\rho_\a=0$. 
Thus we have 
$$
e^{-\beta_\R}\frac{d\rho_\a}{\rho_\a}e^{\beta_\R}=\frac{d\rho_\a}{\rho_\a}-[V_1\w V_2,\, \frac{d\rho_\a}{\rho_\a}]=\frac{d\rho_\a}{\rho_\a}.
$$
Thus we have 
\bgn{align*}
d\J_{\beta_\R}d\log\rho_\a\cdot\ol\psi_{\beta_\R}=&de^{\beta_\R}\J_Je^{-\beta_\R}\frac{d\rho_\a}{\rho_\a}e^{\beta_\R}\cdot\ol\psi\\
=&de^{\beta_\R}\J_J \frac{d\rho_\a}{\rho_\a}\cdot\ol\psi =de^{\beta_\R}\J_J \frac{d\rho_\a}{\rho_\a}e^{-\beta_\R}\ol\psi_{\beta_\R}\\
=&d\J_J\frac{d\rho_\a}{\rho_\a}\cdot\ol\psi_{\beta_\R}+d [\beta_\R,\, \J_J\frac{d\rho_\a}{\rho_\a}]\cdot\ol\psi_{\beta_\R}\\
=&d J(\frac{d\rho_\a}{\rho_\a})\w\psi_{\beta_\R}+d( V_1\lan V_2, \, J\frac{d\rho_\a}{\rho_\a}\ran- V_2\lan V_1,\,
J\frac{d\rho_\a}{\rho_\a}\ran)\cdot\ol\psi_{\beta_\R}\\
=&dJ\frac{d\rho_\a}{\rho_\a}\cdot\psi_{\beta_\R}+\sqrt{-1}d\(\lan V_2, J\frac{d\rho_\a}{\rho_\a}\ran i_{V_1}\ome 
-\lan V_1, J \frac{d\rho_\a}{\rho_\a}\ran i_{V_2}\ome\)\cdot\ol\psi_{\beta_\R}
\end{align*}
Thus we have 
\bgn{align}
&-2d\J_{\beta_\R}\eta_\a\cdot\ol\psi_{\beta_\R}+d\J_{\beta_\R}d\log\rho_\a\cdot\ol\psi_{\beta_\R}\\&=dJ\frac{d\rho_\a}{\rho_\a}\cdot\ol\psi_{\beta_\R}
+\sqrt{-1}d\(\lan V_2, J\frac{d\rho_\a}{\rho_\a}\ran i_{V_1}\ome 
-\lan V_1, J \frac{d\rho_\a}{\rho_\a}\ran i_{V_2}\ome\)\cdot\ol\psi_{\beta_\R}\notag
\end{align}
As in Definition \ref{def:GRic and GR}, 
$-2d\J_{\beta_\R}\eta_\a\cdot\ol\psi_{\beta_\R}+d\J_{\beta_\R}d\log\rho_\a\cdot\ol\psi_{\beta_\R}$
is written as 
$$
-2d\J_{\beta_\R}\eta_\a\cdot\ol\psi_{\beta_\R}+d\J_{\beta_\R}d\log\rho_\a\cdot\ol\psi_{\beta_\R}
=(P-\sqrt{-1}Q)\cdot\ol \psi_{\beta_\R}, 
$$
where $P=GRic$ and $Q$ are real $2$-forms.\\
Since 
$\sqrt{-1}d\(\lan V_2, J\frac{d\rho_\a}{\rho_\a}\ran i_{V_1}\ome 
-\lan V_1, J \frac{d\rho_\a}{\rho_\a}\ran i_{V_2}\ome\)$ is a pure imaginary $2$-form and 
$\dstyle{dJ\frac{d\rho_\a}{\rho_\a}} $ is a real $2$-form, 
we obtain 
$$
\text{\rm GRic} =-dJd\log\rho_\a 
$$
Since $(X,J, \ome)$ is a K\"ahler-Einstein manifold, we also have 
$$-dJ\log\rho_\a 
=\lam \ome.$$ 
Since $\psi_{\beta_\R}=\exp{(-d\mu_{T,1}\w\mu_{T, 2}+\sqrt{-1}\ome)},$
 we have GRic$=\lam \ome$.
\end{proof}
Let $X=(M, J)$ be a compact complex surface with effective anticanonical divisor.
Let $\beta$ be a nontrivial section of $K^{-1}$. Then $\beta$ is a holomorphic Poisson structure.
We denote by  $\J_\beta$ Poisson deformations of generalized complex structures. 
Then from the stability theorem of generalized K\"ahler structures, there is a generalized K\"ahler structure $(\J_\beta , \J_\psi)$, where 
$\psi=e^{B+\sqrt{-1}\ome}$ is a $d$-closed, nondegenerate, pure spinor.
We denote by $D=\{\beta=0\}$  the divisor given by zero of $\beta$.
Then we have 
\bgn{proposition}
Let $\beta$ be a Poisson structure on $X=\Bbb C P^2$ which is an anticanonical divisor $D$ given by three lines in general position. Then there exists a
generalized K\"ahler-Einstein structure $(\J_\beta, \J_\psi)$ such that 
$$
\text{\rm GRic}=3\ome,
$$
where $\psi=e^{b+\sqrt{-1}\ome}$.
\end{proposition}
\bgn{proof}
In our case, Poisson structure $\beta$ is given by an action of $2$-dimensional torus preserving the K\"ahler structure 
of $\C P^2$.  
Then the result follows from Proposition \ref{prop: GKE}.
\end{proof}
\bgn{proposition}
Let $(M,J, \ome)$ be a toric K\"ahler-Einstein manifold of dimension $m$. 
Then there exist deformations of nontrivial generalized K\"ahler-Einstein structures from 
the ordinary K\"ahler-Einstein structure, where $m\geq 2$. 
\end{proposition}
\bgn{proof}
Since $(M, J, \ome)$ is a toric K\"ahler-Einstein manifold, there exists an action of $T^m$ preserving 
the K\"ahler structure.
Then the result follows from Proposition \ref{prop: GKE}
\end{proof}

\medskip
\noindent
E-mail address: goto@math.sci.osaka-u.ac.jp\\
\noindent
Department of Mathematics, Graduate School of Science,\\
\noindent Osaka University Toyonaka, Osaka 560-0043, JAPAN 
\end{document}